\documentclass[12pt]{amsart}

\usepackage{amssymb}
\usepackage{graphicx}
\usepackage{float}

\DeclareMathAlphabet{\mathpzc}{OT1}{pzc}{m}{it}

\newtheorem{theorem}{Theorem}[section]
\newtheorem{lemma}[theorem]{Lemma}

\theoremstyle{definition}

\theoremstyle{remark}

\numberwithin{equation}{section}

\begin{document}

\title{Alexander--equivalent Zariski pairs of irreducible sextics}

\author{Christophe Eyral and Mutsuo Oka}

\address{Department of Mathematical Sciences, University of Aarhus, Building 1530, Ny Munkegade, DK--8000 Aarhus C, Denmark}  
\email{eyralchr@yahoo.com}
\address{Department of Mathematics, Tokyo University of Science, 26 Wakamiya--cho, Shinjuku--ku, Tokyo 162--8601, Japan}   
\email{oka@rs.kagu.tus.ac.jp}

\subjclass[2000]{14H30 (14H20, 14H45, 14H50).}

\keywords{Topology of plane curves; Fundamental group; Zariski pairs; Alexander polynomial.}

\thanks{}  

\begin{abstract}
The existence of Alexander--equivalent Zariski pairs dealing with irreducible curves of degree 6 was proved by A. Degtyarev. However, up~to~now, no explicit example of such a pair was available (only the existence was known). In this paper, we construct the first concrete example.
\end{abstract}

\maketitle

\markboth{C. Eyral and M. Oka}{Alexander--equivalent Zariski pairs}  

\section{Introduction}

Let $\mathcal{M}(\Sigma,d)$ be the moduli space of reduced curves of
degree $d$ in $\mathbb{CP}^2$ with a prescribed configuration of
singularities $\Sigma$.\footnote{By such a moduli space we mean the quotient space $\mathcal{C}(\Sigma,d)/\text{PGL}(3,\mathbb{C})$ of the space $\mathcal{C}(\Sigma,d)$ of reduced plane curves with degree $d$ and set of singularities $\Sigma$ by the `standard' group action of $\text{PGL}(3,\mathbb{C})$.} A pair of curves $(C,C')$ in $\mathcal{M}(\Sigma,d)$ is said to be a \emph{Zariski pair} if it satisfies the following two conditions (cf.~\cite{A}):
\begin{enumerate}
\item \label{c1} $C$ and $C'$ have the same \emph{combinatoric}, that is, there exist regular neighbourhoods $T(C)$ and $T(C')$ of $C$ and $C'$ respectively such that the pairs $(T(C),C)$ and $(T(C'),C')$ are homeomorphic;
\item \label{c2} the pairs $(\mathbb{CP}^2,C)$ and $(\mathbb{CP}^2,C')$ are not homeomorphic.
\end{enumerate}
It is easy to check that if both $C$ and $C'$ are \emph{irreducible} then the first condition is always satisfied. The first Zariski pair appears in the works by O. Zariski \cite{Z1,Z2,Z3} (see also \cite{A} and \cite{O3}). The members of the pair are irreducible curves of degree~$6$, which is the smallest degree for which Zariski pairs exist.

The existence of a Zariski pair in a moduli space gives an information about its connected components. Indeed, if $\mathcal{M}(\Sigma,d)$ has a Zariski pair $(C,C')$, then $C$ and $C'$ necessarily belong to different connected components (cf.~\cite{Z4,Z5,LR}) --- in particular, if $\mathcal{M}(\Sigma,d)$ has a Zariski pair, then it is not connected. (The converse statement is not clear. Two curves coming from two different connected components may have the same embedded topology.)

Now to check whether a given a pair of curves $C,C' \in \mathcal{M}(\Sigma,d)$ with the same combinatoric is a Zariski pair, one can first try to calculate the generic Alexander polynomials of the curves. If these polynomials are different, then the curves do not have the same embedded topology, and $(C,C')$ is a Zariski pair. However, it may happen that these polynomials are the same although $(C,C')$ is a Zariski pair. In this case the pair is said to be \emph{Alexander--equivalent}. The first example of such a pair was given by E. Artal Bartolo and J. Carmona Ruber~\cite{AC} for \emph{reducible} curves of degree 7. The first examples dealing with \emph{irreducible} curves are due to the second author \cite{O1} (curves of degree $12$) and \cite{O2} (curves of degree $8$). 
In~\cite{D1}, A. Degtyarev proved that Alexander--equivalent Zariski pairs also appear on irreducible curves of degree $6$. However he did not give any explicit example (only the existence is proved). The aim of the present paper is to construct a concrete example of such a pair. 

\section{Statement of the result}

\begin{figure}[H]
\includegraphics[width=8cm,height=5cm]{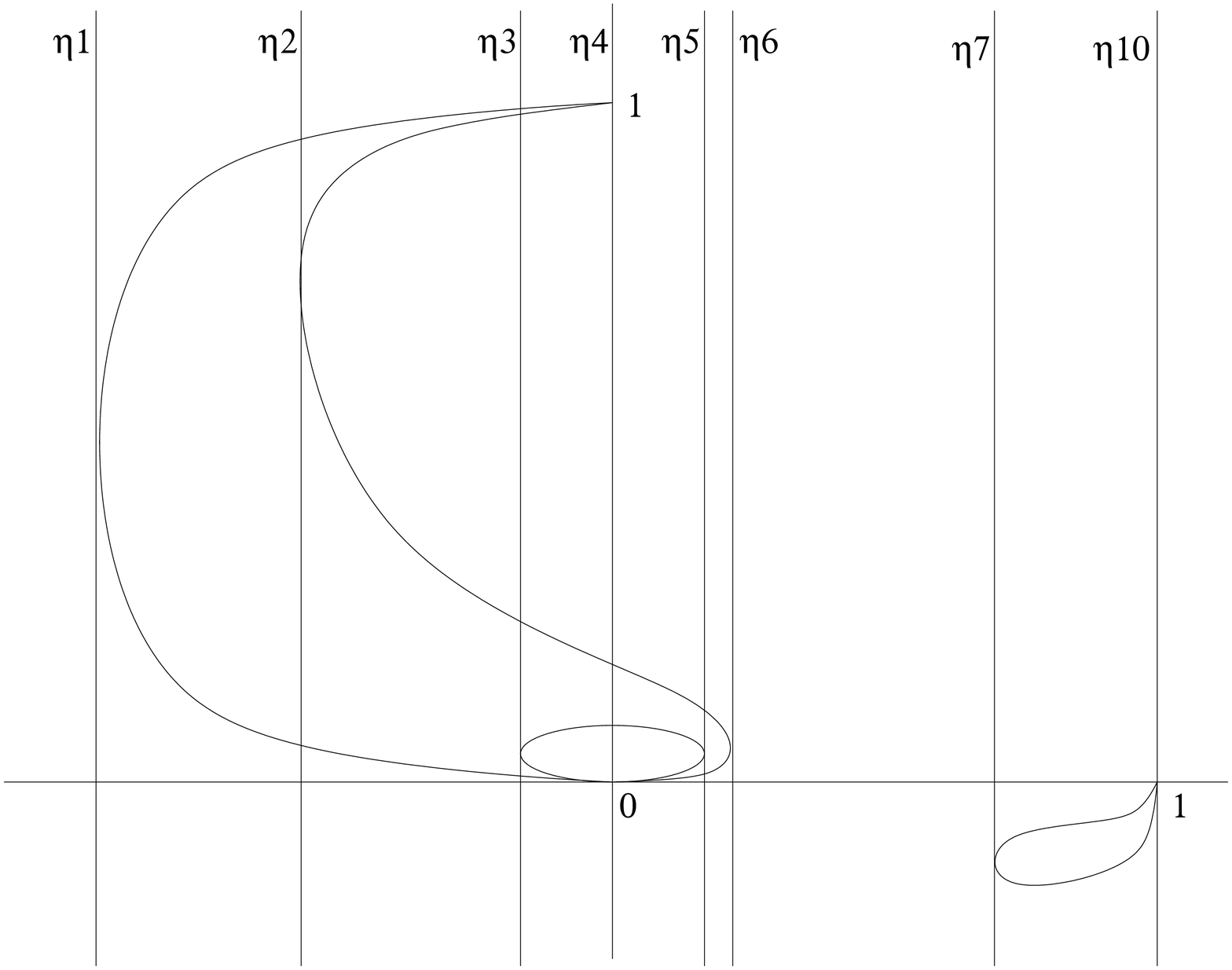}
\caption{\label{figa92a4rps}$\{(x,y)\in\mathbb{R}^2\ ;\ f(x,y)=0\}$}
\end{figure}

\begin{figure}[H]
\includegraphics[width=8cm,height=5cm]{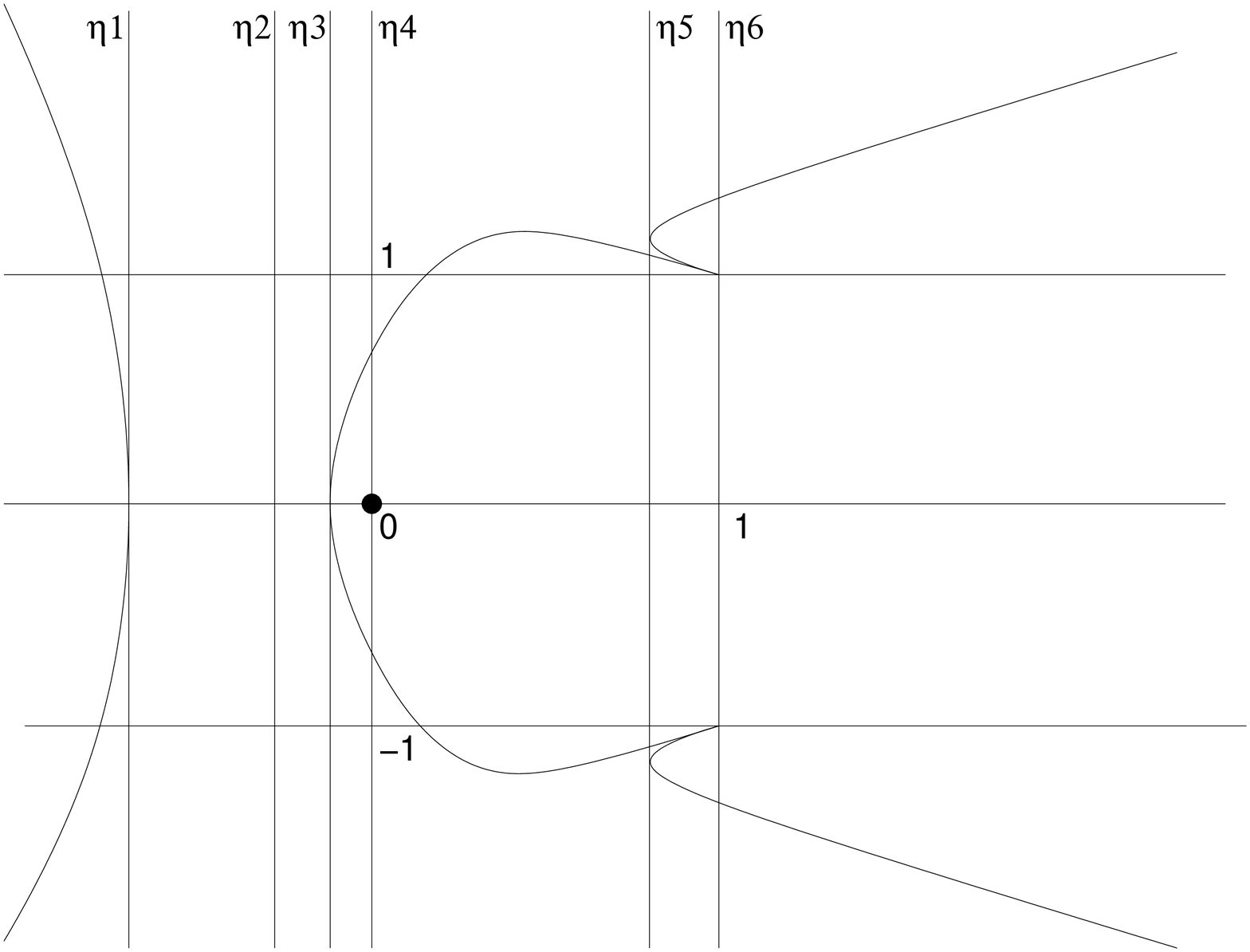}
\caption{\label{figa92a4narps}$\{(x,y)\in\mathbb{R}^2\ ;\ f'(x,y)=0\}$}
\end{figure}

Let $(X\colon Y\colon Z)$ be homogeneous coordinates on $\mathbb{CP}^2$ and $(x,y)$ the affine coordinates defined by $x:=X/Z$ and $y:=Y/Z$ on $\mathbb{CP}^2 \setminus \{Z=0\}$, as usual. We consider the projective curves $C$ and $C'$ in $\mathbb{CP}^2$ defined by the affine equations $f(x,y)=0$ and $f'(x,y)=0$ respectively, where 
\begin{eqnarray*}
f(x,y) & := & {\frac {369}{364}}\,{y}^{6}+{y}^{5}x-{\frac {197}{91}}\,{y}^{5}+{
\frac {207}{182}}\,{y}^{4}{x}^{2}-{\frac {185}{91}}\,{y}^{4}x+{\frac {
235}{182}}\,{y}^{4}+\\
& & {\frac {87}{182}}\,{y}^{3}{x}^{3}-{\frac {255}{182
}}\,{y}^{3}{x}^{2}+{\frac {97}{91}}\,{y}^{3}x-{\frac {1}{7}}\,{y}^{3}+{\frac {101
}{364}}\,{y}^{2}{x}^{4}-{\frac {47}{91}}\,{y}^{2}{x}^{3}+\\
& & {\frac {7}{26
}}\,{y}^{2}{x}^{2}-{\frac {3}{91}}\,{y}^{2}x+{\frac {1}{364}}\,{y}^{2}
+{\frac {5}{182}}\,y{x}^{5}-{\frac {11}{182}}\,y{x}^{4}+{\frac {1}{26}}\,y{x}^{3}-\\
& & {\frac {1}{182}}\,y{x}^{2}+{\frac {1}{364}}\,{x}^{6}-{\frac {1}{182}}
\,{x}^{5}+{\frac {1}{364}}\,{x}^{4},\\
&&\\
f'(x,y) & := & -{\frac {4}{3}}\,{y}^{6}+ \left( -{\frac {8}{9}} \, {x}^{2}+4\,x+1 \right) {y}^{4}+ \left( -{\frac {4}{9}}\,{x}^{4}+{\frac {26}{9}} \,{x}^{3}-{\frac {14}{3}}\,{x}^{2}-2\,x \right) {y}^{2}\\
& & +{\frac {1}{9}}\,{x}^{6}+{\frac {2}{9}}\,{x}^{5}-{\frac {17}{9}}\,{x}^{4}+2
\,{x}^{3}+{x}^{2}.
\end{eqnarray*}
Both $C$ and $C'$ are \emph{irreducible sextics} with the set of
singularities $\textbf{A}_9\oplus 2\textbf{A}_4$. (We recall that a point $P$ in a curve $D$ is said to be an $\textbf{A}_n$--singularity ($n\geq 1$) if the germ $(D,P)$ is topologically equivalent to the germ at the origin of the curve defined by $x^2+y^{n+1}=0$.)
For the curve $C$, the $\textbf{A}_4$--singularities are located at $(1,0)$ and $(0,1)$ while the $\textbf{A}_9$--singu\-lari\-ty is at $(0,0)$.  
For $C'$, the $\textbf{A}_4$--singularities are at $(1,\pm 1)$
and the $\textbf{A}_9$--sing\-ularity is at $(0,0)$.
The \emph{real} plane sections of $C$ and $C'$ are shown in
Figure~\ref{figa92a4rps} and \ref{figa92a4narps} respectively. (In the figures we do not respect the numerical scale.) The curve $C'$ is symmetric with respect to the $x$--axis. The curve $C$ has no particular symmetry.
Notice that after the 
analytic change of coordinates 
\begin{equation*}
(x,y) \mapsto \left(x-\frac{1}{3}\, y^4+y^2,y\right),
\end{equation*}
the equation of $C'$ near the origin takes the form
\begin{equation*}
x^2+\frac{4}{27}\, y^{10}+\mbox{higher terms} = 0.
\end{equation*}
As the leading term $x^2+(4/27)\, y^{10}$ is positive on $\mathbb{R}^2\setminus \{(0,0)\}$, the origin is an isolated point of the real plane section of $C'$ (cf.~Figure~\ref{figa92a4narps}).
Finally, since $C$ and $C'$ are not of torus type (cf.~\cite{OP} and \cite{P}), their generic Alexander polynomials are trivial (cf.~\cite{D1}).

\begin{theorem}\label{maintheorem}
The fundamental group $\pi_1(\mathbb{CP}^2\setminus C)$ is isomorphic to $\mathbb{Z}/6\mathbb{Z}$ while $\pi_1(\mathbb{CP}^2\setminus C')$ is isomorphic to $\mathbb{D}_{10}\times (\mathbb{Z}/3\mathbb{Z})$, where $\mathbb{D}_{10}$ is the dihedral group of order $10$. In particular $(C,C')$ is an Alexander--equivalent Zariski pair and the moduli space of irreducible sextics with the set of singularities $\textbf{\emph{A}}_9\oplus 2\textbf{\emph{A}}_4$ has at least two connected components.\footnote{Though the existence of the structure of an algebraic variety on a moduli space is not always obvious, the moduli space we consider here has such a structure. The last assertion in the theorem then implies this moduli space has at least two irreducible components as well.}
\end{theorem}

Theorem \ref{maintheorem} is proved in sections \ref{abelian} and \ref{nonabelian} below.

The curve $C'$ is an example of so--called \emph{$\mathbb{D}_{10}$--sextics}. (A $\mathbb{D}_{10}$--sextic is a non--torus irreducible sextic with simple singularities and whose fundamental group\footnote{We always mean the fundamental group of the complement of the curve.} factors to the dihedral group $\mathbb{D}_{10}$.) The existence of $\mathbb{D}_{10}$--sextics was first proved, purely arithmetically, by A. Degtyarev \cite{D1} who showed that there exist exactly 8 equisingular deformation families of such curves, one family for each of the following sets of singularities:
\begin{center}
$4\textbf{A}_4,\ 4\textbf{A}_4\oplus\textbf{A}_1,\
4\textbf{A}_4\oplus2\textbf{A}_1,\ 4\textbf{A}_4\oplus\textbf{A}_2,$
$\textbf{A}_9\oplus2\textbf{A}_4,\
\textbf{A}_9\oplus2\textbf{A}_4\oplus\textbf{A}_1,\
\textbf{A}_9\oplus2\textbf{A}_4\oplus\textbf{A}_2,\
2\textbf{A}_9.$
\end{center}
First explicit examples and fundamental groups of $\mathbb{D}_{10}$--sextics were given in \cite{D2} (see also \cite{EO} for the sets of singularities $4\textbf{A}_4$ and $4\textbf{A}_4\oplus\textbf{A}_1$). 

Furthermore, Degtyarev also observed in \cite{D1} that the configurations 
\begin{displaymath}
4\textbf{A}_4,\ 4\textbf{A}_4\oplus\textbf{A}_1,\
\textbf{A}_9\oplus2\textbf{A}_4,\
\textbf{A}_9\oplus2\textbf{A}_4\oplus\textbf{A}_1,\
2\textbf{A}_9 
\end{displaymath}
can be realized not only by $\mathbb{D}_{10}$--sextics but also by irreducible non--torus sextics $D$ for which the group $\pi_1(\mathbb{CP}^2\setminus D)$ does \emph{not} factors to~$\mathbb{D}_{10}$. (In particular, since the Alexander polynomial of an irreducible non--torus sextic with simple singularities is always trivial, these 5 sets of singularities give rise to Alexander--equivalent Zariski pairs of irreducible sextics.) However, Degtyarev did not give any explicit equation of such a curve $D$ (only the existence is proved) and did not calculate the group $\pi_1(\mathbb{CP}^2\setminus D)$. The curve $C$ in Theorem \ref{maintheorem} is the first explicit example of such a curve --- i.e., a curve whose fundamental group does not admit a dihedral quotient although its set of singularities can be realized by a $\mathbb{D}_{10}$--sextic as well. (In particular, the pair $(C,C')$ in Theorem \ref{maintheorem} is also the first concrete example of an Alexander--equivalent Zariski pair dealing with irreducible sextics.)

\section{Fundamental group of $\mathbb{CP}^2\setminus C$}\label{abelian}

In this section, we show that $\pi_1(\mathbb{CP}^2\setminus C)\simeq \mathbb{Z}/6\mathbb{Z}$. In fact, it suffices to prove that $\pi_1(\mathbb{CP}^2\setminus C)$ is abelian. Indeed, by Hurewicz's theorem, if $\pi_1(\mathbb{CP}^2\setminus C)$ is abelian, then it is isomorphic to first integral homology group $H_1(\mathbb{CP}^2\setminus C)$ and it is well known that $H_1(\mathbb{CP}^2\setminus C)\simeq \mathbb{Z}/6\mathbb{Z}$.

To show that $\pi_1(\mathbb{CP}^2\setminus C)$ is abelian, we use Zariski--van Kampen's theorem with the pencil given by the vertical lines $L_{\eta}\colon x=\eta$, $\eta\in\mathbb{C}$ (cf.~\cite{Z1} and \cite{vK}). We always take the point $(0\colon 1\colon 0)$ as base point for our fundamental groups. This point is nothing but the axis of the pencil, which is also the point at infinity of the lines $L_{\eta}$. Note that it does not belong to the curve.

The discriminant $\Delta_y(f)$ of $f$ as a polynomial in $y$ is the  polynomial in $x$ given by
\begin{eqnarray*}
\Delta_y(f)(x) =  a_0\, {x}^{15}\, \, (x-1)^{7}\, (858898351\,{x
}^{8}-1278576626\,{x}^{7}-359900737\,{x}^{6}+\\
1017975356\,{x}^{5}-56181608\,{x}^{4}- 170653568\,{x}^{3}+ 2388080\,{x}^{2}+2072000\,x-96000),
\end{eqnarray*}
where $a_0\in\mathbb{Q}\setminus\{0\}$.
This polynomial has exactly 10 distinct complex roots:\smallskip

$\eta_1 \approx  - 0.7408$,
\quad $\eta_2 \approx - 0.3914$,
\quad $\eta_3 \approx - 0.1309$,
\quad $\eta_4 = 0$,

$\eta_5 \approx  0.0598$,
\quad $\eta_6 \approx  0.0778$,
\quad $\eta_7 \approx  0.6274$,

$\eta_8 \approx  0.9933 - i\, 0.1446$,
\quad $\eta_9 =\bar\eta_8\approx 0.9933 + i\, 0.1446$,
\quad $\eta_{10} = 1$.

\smallskip

\noindent
The singular lines of the pencil are the lines
$L_{\eta_j}$ ($1\leq j\leq 10$) corresponding to these 10 roots. The
lines $L_{\eta_4}$ and $L_{\eta_{10}}$ intersect the curve at its singular points. All the other singular lines are tangent to $C$. See Figure~\ref{figa92a4rps}.

\begin{figure}[H]
\includegraphics[width=8cm,height=5cm]{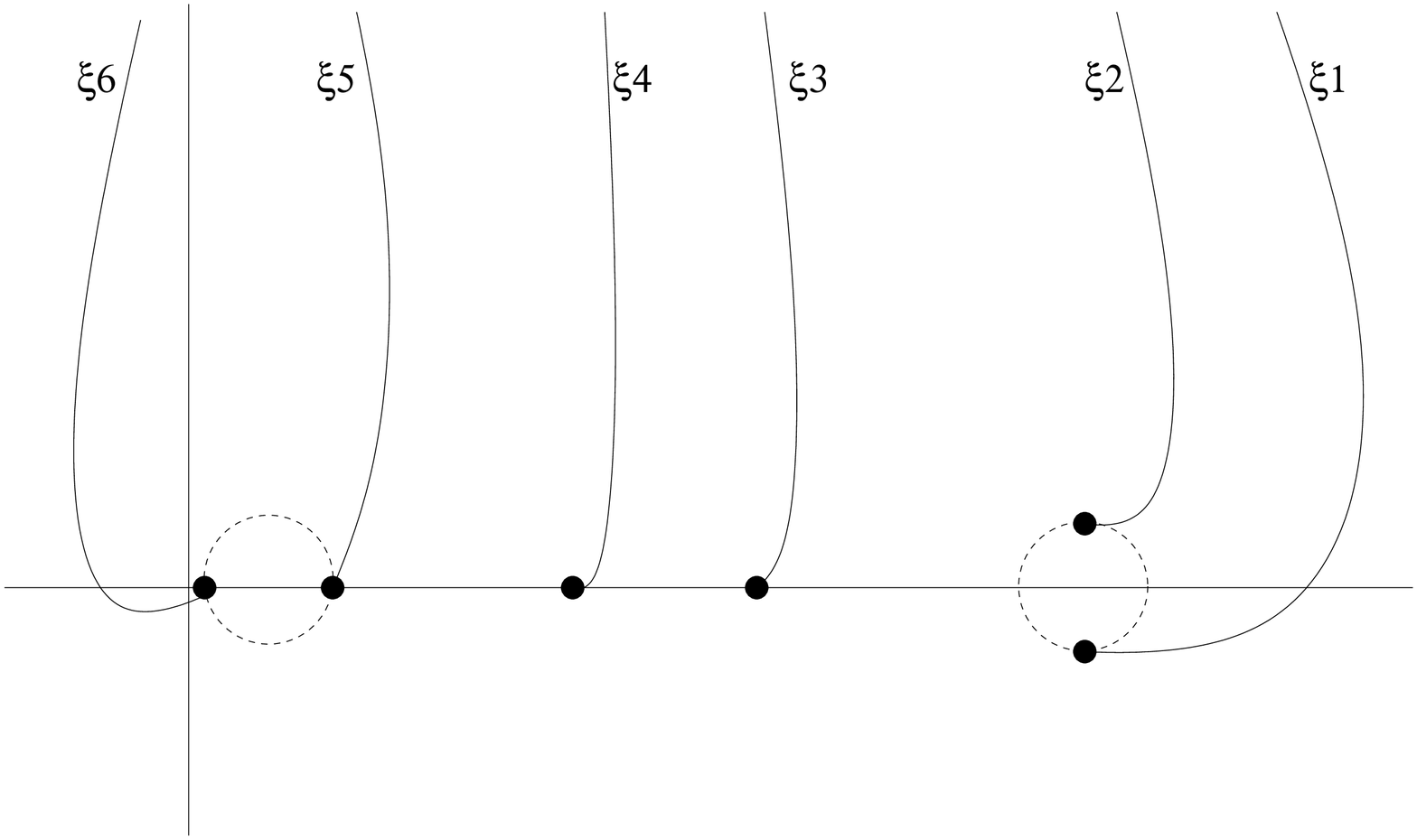}
\caption{\label{figa92a4eta4}Generators at $x=\eta_4+\varepsilon$}
\end{figure}

We consider the generic line $L_{\eta_4+\varepsilon}$ and
choose generators $\xi_1,\ldots,\xi_6$ of the fundamental group 
$\pi_1(L_{\eta_4+\varepsilon}\setminus C)$ as in Figure~\ref{figa92a4eta4}, where $\varepsilon>0$ is small enough. The $\xi_k$'s ($1\leq k\leq 6$) are \emph{lassos oriented counter--clockwise} around the 6 intersection points of the line $L_{\eta_4+\varepsilon}$ with the curve --- i.e., the 6 complex roots of the equation $f(\eta_4+\varepsilon,y)=0$. (In the figures, a lasso is represented by a path ending with a bullet.)
The~Zariski--van Kampen theorem says that
\begin{displaymath}
\pi_1(\mathbb{CP}^2\setminus C) \simeq
\pi_1(L_{\eta_4+\varepsilon}\setminus C) \big/ G,
\end{displaymath}
where $G$ is the normal subgroup of $\pi_1(L_{\eta_4+\varepsilon}
\setminus C)$ generated by the monodromy relations associated with the
singular lines of the pencil. To find these relations, we
fix a `standard' system of generators $\sigma_1,\ldots,\sigma_{10}$ for
the fundamental group $\pi_1(\mathbb{C}\setminus \{\eta_1,
\ldots, \eta_{10}\})$ as follows. Each $\sigma_j$ is a lasso (oriented counter--clockwise) around $\eta_j$ with base point
$\eta_4+\varepsilon$. For $j\not= 8,9$, the tail of $\sigma_j$ is a union of real segments and half--circles around the exceptional parameters $\eta_l$ ($l\not=j$) located in the real axis between the base point $\eta_4+\varepsilon$ and $\eta_j$. Its head is the circle $\mathbb{S}_\varepsilon(\eta_j)$ with centre $\eta_j$ and radius $\varepsilon$. The lasso $\sigma_8$ corresponding to the non--real root $\eta_8$ is given by $\zeta\theta\zeta^{-1}$, where $\theta$ is the loop obtained by moving $x$ once on the circle $\mathbb{S}_\varepsilon(\eta_8)$, starting at $\Re(\eta_8)+i\, (\Im(\eta_8)+\varepsilon)$, while $\zeta$ is the path obtained when $x$ moves on the real axis from $\eta_4+\varepsilon$ to $\eta_5-\varepsilon$, makes half--turn on the circle $\mathbb{S}_\varepsilon(\eta_5)$, from $\eta_5-\varepsilon$ to $\eta_5+\varepsilon$, moves on the real axis from $\eta_5+\varepsilon$ to $\eta_6-\varepsilon$, makes half--turn on the circle $\mathbb{S}_\varepsilon(\eta_6)$, from $\eta_6-\varepsilon$ to $\eta_6+\varepsilon$, moves on the real axis from $\eta_6+\varepsilon$ to $\eta_7-\varepsilon$, makes half--turn on the circle $\mathbb{S}_\varepsilon(\eta_7)$, from $\eta_7-\varepsilon$ to $\eta_7+\varepsilon$, moves on the real axis from $\eta_7+\varepsilon$ to $\Re(\eta_8)$, and finally moves in a straight line from $\Re(\eta_8)$ to $\Re(\eta_8)+i\, (\Im(\eta_8)+\varepsilon)$. (Here $\Re(\eta_8)$ and $\Im(\eta_8)$ denote the real and the imaginary parts of $\eta_8$ respectively.) The lasso $\sigma_9$ is defined similarly from a loop $\theta$ and a path $\zeta$ meeting at $\Re(\eta_9)+i\, (\Im(\eta_9)-\varepsilon)$.
The monodromy relations around the singular line $L_{\eta_j}$ are
obtained by moving  the generic `fibre' $F\simeq L_{\eta_4+\varepsilon}
\setminus C$ isotopically `above' the loop $\sigma_j$, and by identifying each $\xi_k$ ($1\leq k\leq 6$) with its image by the terminal homeomorphism of this isotopy. (For details see \cite{Z1} and~\cite{vK}.) 

The remaining of the proof is to determine these relations.

\begin{figure}[H]
\includegraphics[width=8cm,height=5cm]{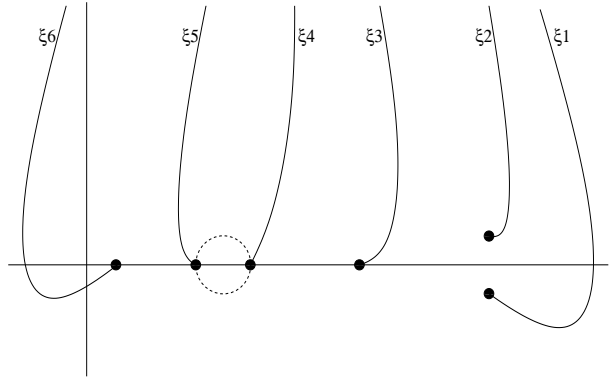}
\caption{\label{figa92a4eta5}Generators at $x=\eta_5-\varepsilon$}
\end{figure}

\subsection*{Monodromy relations at $x=\eta_5$}
When $x$ moves on the real axis from $\eta_4+\varepsilon$ to $\eta_5-\varepsilon$, the 6 complex roots of the equation (in $y$) $f(x,y)=0$ (and, consequently, the 6 generators $\xi_1,\ldots\xi_6$) are deformed as in Figure~\ref{figa92a4eta5}. The singular line $L_{\eta_5}$ is tangent to the curve at the simple point $P\approx(\eta_5,0.0095)$, and the intersection multiplicity $I(L_{\eta_5},C;P)$ of this line with the curve at this point is 2. Therefore, by the implicit functions theorem, the germ $(C,P)$  is given by 
\begin{displaymath}
x-\eta_5=b_0\, (y-0.0095)^2+\mbox{higher terms},
\end{displaymath}
where $b_0\not=0$. It follows that when $x$ runs once counter--clockwise on the circle $\mathbb{S}_\varepsilon(\eta_5)$, starting at $\eta_5-\varepsilon$, the variable $y$ makes half--turn counter--clockwise on the dotted circle around $0.0095$ (cf.~Figure~\ref{figa92a4eta5}). The monodromy relation at $x=\eta_5$ is then given by 
\begin{equation}\label{rela92a4eta5}
\xi_5=\xi_4.
\end{equation}

\subsection*{Monodromy relations at $x=\eta_6$}
\begin{figure}[H]
\includegraphics[width=8cm,height=5cm]{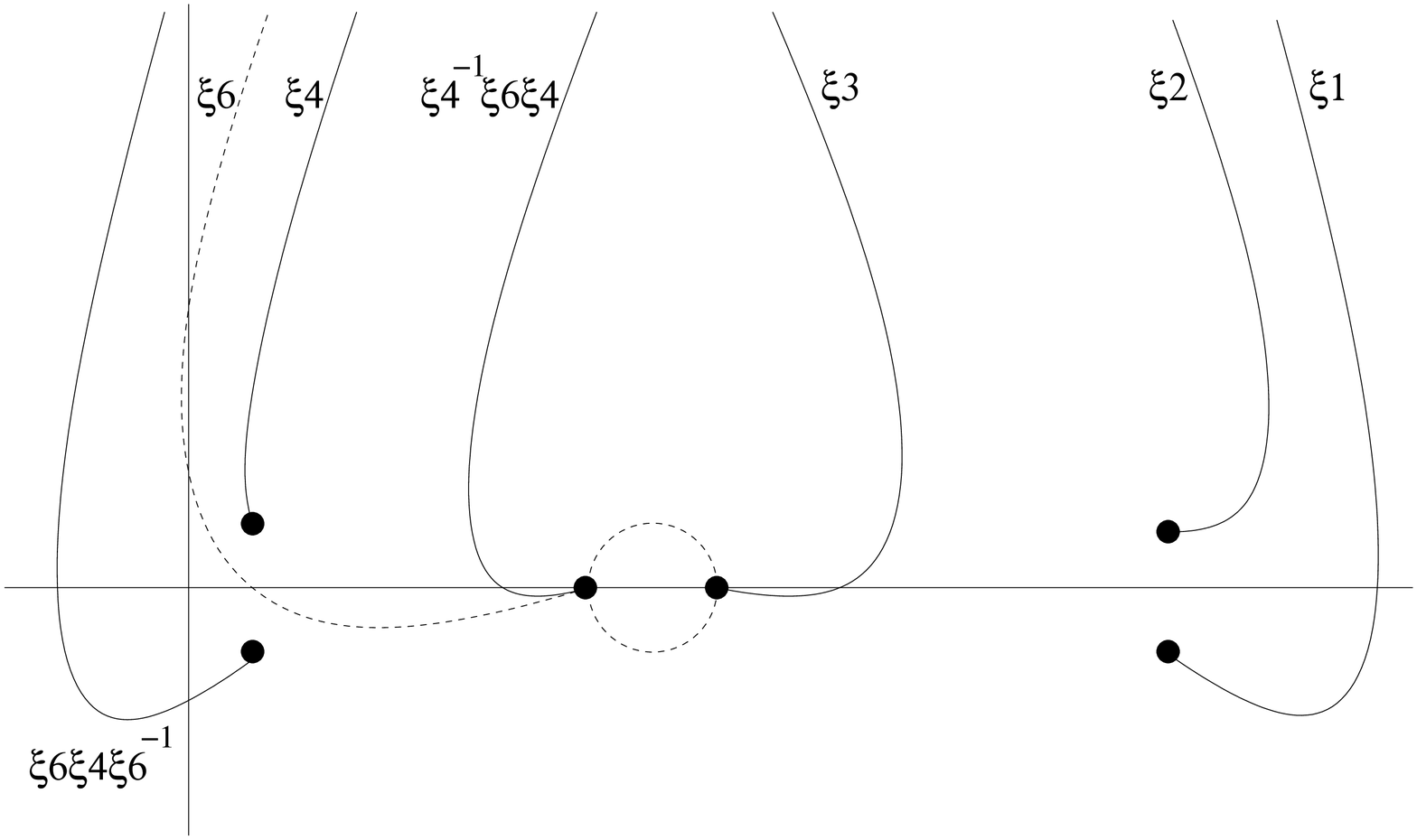}
\caption{\label{figa92a4eta6}Generators at $x=\eta_6-\varepsilon$}
\end{figure}

At $x=\eta_5-\varepsilon$, the generators are as in Figure~\ref{figa92a4eta5}.
In Figure~\ref{figa92a4eta6}, we show how the $\xi_k$'s are deformed when $x$ first makes half--turn counter--clockwise on the circle $\mathbb{S}_\varepsilon(\eta_5)$, from $\eta_5-\varepsilon$ to $\eta_5+\varepsilon$, and then moves on the real axis from $\eta_5+\varepsilon$ to $\eta_6-\varepsilon$. The singular line $L_{\eta_6}$ is also tangent to $C$ at a simple point $P'$ and $I(L_{\eta_6},C;P')=2$. Therefore, as above, the monodromy relation we are looking for is simply given by 
\begin{equation}\label{rela92a4eta6}
\xi_3=\xi_4^{-1}\xi_6\xi_4.
\end{equation}

\subsection*{Monodromy relations at $x=\eta_4$}
The singular line $L_{\eta_4}$ passes through the singular points $(0,0)$ and $(0,1)$ which are singularities of type $\textbf{A}_9$ and $\textbf{A}_4$ repectively. At $(0,0)$, the curve has two branches $K_-$ and $K_+$ given by
\begin{equation*}
K_\pm:\quad y=x^2+5\, x^3+51\, x^4+(503\pm 32 \sqrt{6})\, x^5 +\mbox{higher terms.}
\end{equation*}
It follows that when $x$ runs once counter--clockwise on the circle $\mathbb{S}_\varepsilon (\eta_4)$, starting at $\eta_4+\varepsilon$, the generators $\xi_6$ and $\xi_5=\xi_4$ make $5$ turns counter--clockwise on the corresponding dotted circle (cf.~Figure~\ref{figa92a4eta4}).
The monodromy relation around $L_{\eta_4}$ that comes from the singular point $(0,0)$ is then given by
\begin{equation}\label{rela92a4eta400}
\xi_4=(\xi_6\xi_4)^5\cdot\xi_4\cdot(\xi_6\xi_4)^{-5}.
\end{equation}
At $(0,1)$, a Puiseux parametrization of $C$ is given by
\begin{equation*}
x=t^2,\quad y=1-\frac{1}{2}\, t^4+\frac{1}{10}\,i\, \sqrt{5}\, t^5+\mbox{higher terms.}
\end{equation*}
Hence, when $x$ goes once counter--clockwise on the circle $\mathbb{S}_\varepsilon (\eta_4)$, starting at $\eta_4+\varepsilon$, the generators $\xi_1$ and $\xi_2$ make $(5/2)$--turn counter--clockwise on the corresponding dotted circle (cf.~Figure~\ref{figa92a4eta4}). The monodromy relation around $L_{\eta_4}$ that comes from the singular point $(0,1)$ is then given by 
\begin{equation}\label{rela92a4eta401}
\xi_1=(\xi_2\xi_1)^2\cdot\xi_2\cdot(\xi_2\xi_1)^{-2}.
\end{equation}

\subsection*{Monodromy relations at $x=\eta_3$}
\begin{figure}[H]
\includegraphics[width=8cm,height=5cm]{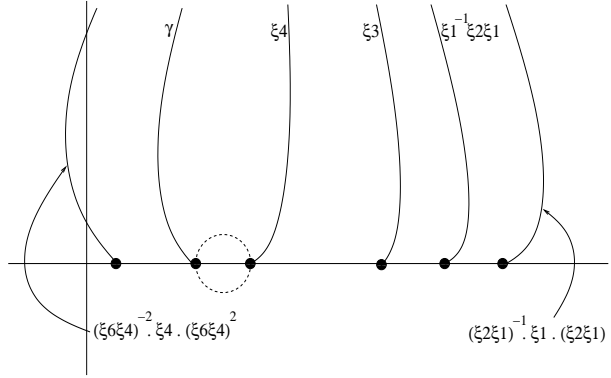}
\caption{\label{figa92a4eta3}Generators at $x=\eta_3+\varepsilon$}
\end{figure}

At $x=\eta_4+\varepsilon$, the generators are as in Figure~\ref{figa92a4eta4}.
Now, when $x$ makes half--turn counter--clockwise on the circle $\mathbb{S}_\varepsilon (\eta_4)$, from $\eta_4+\varepsilon$ to from $\eta_4-\varepsilon$, and then moves on the real axis from $\eta_4-\varepsilon$ to $\eta_3+\varepsilon$, the $\xi_k$'s are deformed as 
in Figure~\ref{figa92a4eta3}, where 
\begin{displaymath}
\gamma:=(\xi_6\xi_4)^{-2}\cdot\xi_4^{-1}\xi_6\xi_4\cdot(\xi_6\xi_4)^2.
\end{displaymath}
The singular line $L_{\eta_3}$ is tangent to the curve at a simple point, with intersection multiplicity 2, and the monodromy relation we are looking for is given by 
\begin{equation}\label{rela92a4eta3}
\xi_4=\gamma.
\end{equation}

\begin{figure}[H]
\includegraphics[width=8cm,height=5cm]{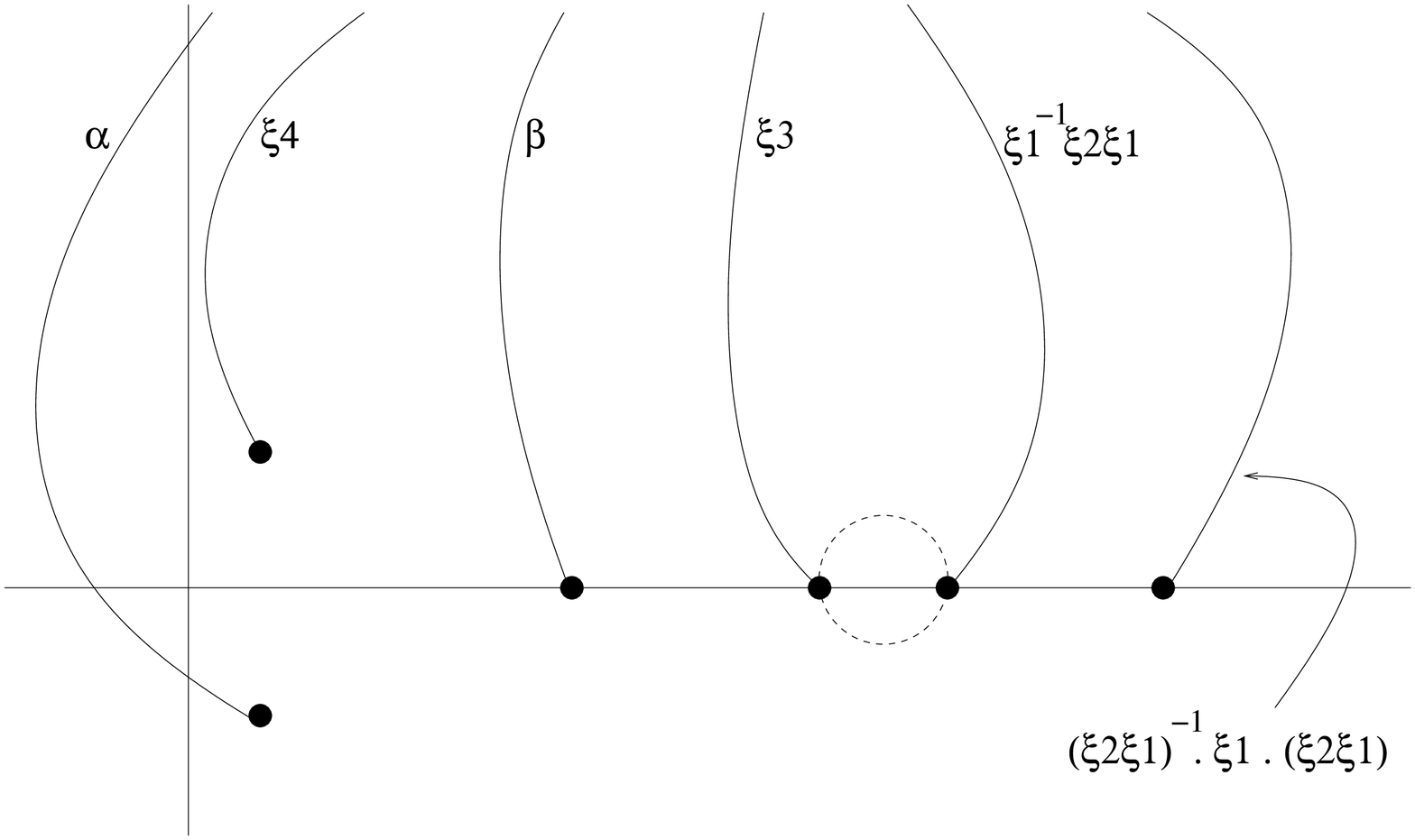}
\caption{\label{figa92a4eta2}Generators at $x=\eta_2+\varepsilon$}
\end{figure}

\subsection*{Monodromy relations at $x=\eta_2$}
When $x$ makes half--turn counter--clockwise on the circle $\mathbb{S}_\varepsilon(\eta_3)$, from $\eta_3+\varepsilon$ to $\eta_3-\varepsilon$, then moves on the real axis from $\eta_3-\varepsilon$ to $\eta_2+\varepsilon$, the $\xi_k$'s are deformed as in
Figure~\ref{figa92a4eta2}, where 
\begin{eqnarray*}
& & \alpha:=\bigl((\xi_6\xi_4)^{-2}\cdot\xi_4\cdot(\xi_6\xi_4)^2\bigr)
\cdot\xi_4\cdot \bigl((\xi_6\xi_4)^{-2} \cdot\xi_4 \cdot (\xi_6\xi_4)^2\bigr)^{-1},\\
& & \beta:=\xi_4^{-1} \cdot (\xi_6\xi_4)^{-2}\cdot \xi_4 \cdot(\xi_6\xi_4)^2 \cdot \xi_4.
\end{eqnarray*}
The monodromy relation at $x=\eta_2$ is also an usual multiplicity 2 tangent relation:
\begin{equation}\label{rela92a4eta2}
\xi_3=\xi_1^{-1}\xi_2\xi_1.
\end{equation}

\begin{figure}[H]
\includegraphics[width=8cm,height=5cm]{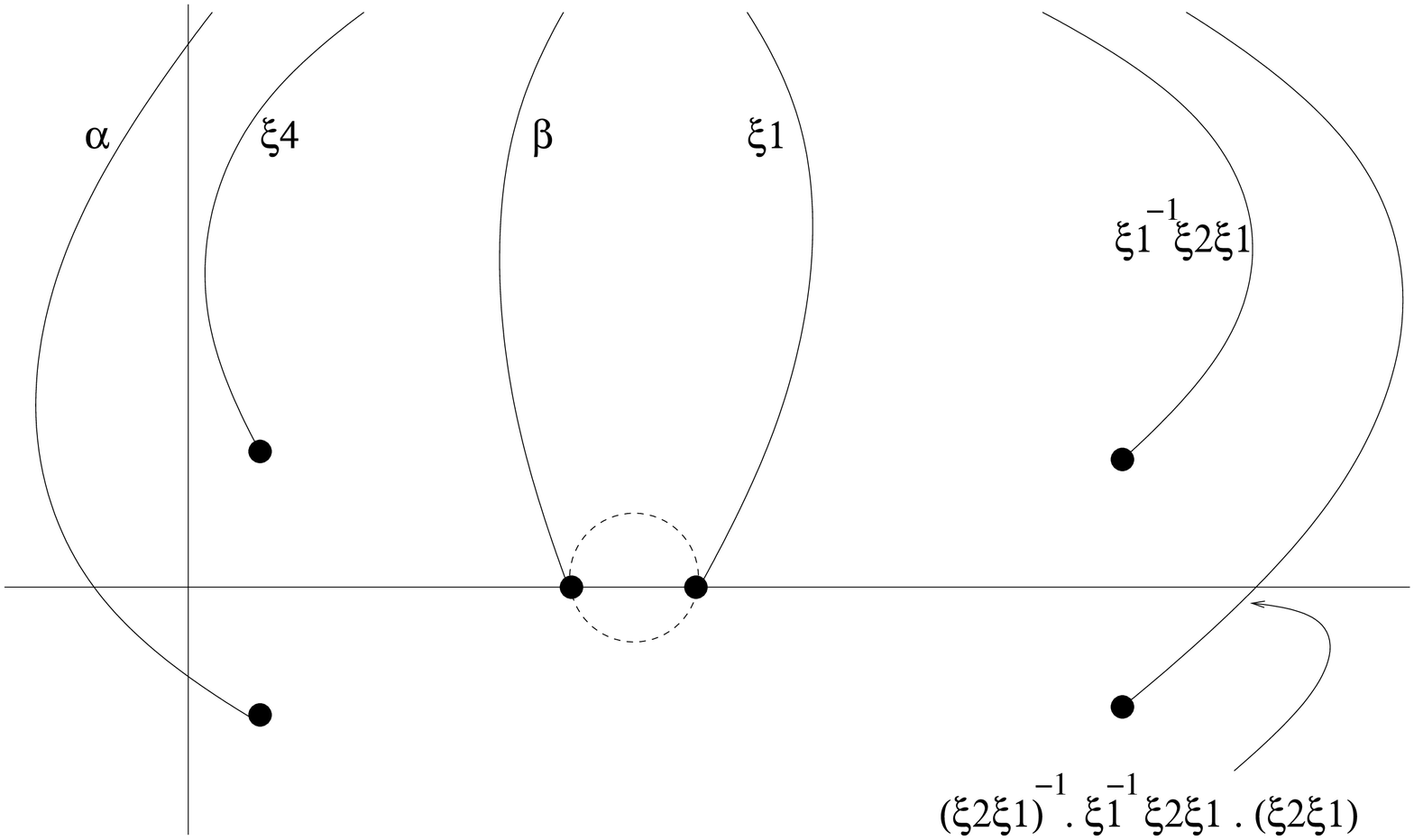}
\caption{\label{figa92a4eta1}Generators at $x=\eta_1+\varepsilon$}
\end{figure}

\subsection*{Monodromy relations at $x=\eta_1$}
In Figure~\ref{figa92a4eta1}, we show how the $\xi_k$'s are deformed when $x$ makes half--turn counter--clockwise on the circle $\mathbb{S}_\varepsilon(\eta_2)$, from $\eta_2+\varepsilon$ to $\eta_2-\varepsilon$, then moves on the real axis from $\eta_2-\varepsilon$ to $\eta_1+\varepsilon$.
The monodromy relation around $L_{\eta_1}$ is a multiplicity 2 tangent relation given by
\begin{equation}\label{rela92a4eta1}
\xi_1=\beta.
\end{equation}
Equivalently $(\xi_6\xi_4)^2 \cdot \xi_4 = \xi_4^{-1}(\xi_6\xi_4)^2 \cdot \xi_4 \xi_1$. Since $(\xi_6\xi_4)^2 \cdot \xi_4 = \xi_4^{-1}(\xi_6\xi_4)^2\cdot\xi_6\xi_4$ (by (\ref{rela92a4eta3})), it follows that
\begin{equation}\label{rela92a4s1}
\xi_4\xi_1=\xi_6\xi_4.
\end{equation}
Combined with (\ref{rela92a4eta6}), this gives
\begin{equation}\label{rela92a4s2}
\xi_3=\xi_1.
\end{equation}
Combined with (\ref{rela92a4eta2}), this in turn implies
\begin{equation}\label{rela92a4s3}
\xi_2=\xi_1.
\end{equation}

\subsection*{Monodromy relations at $x=\eta_7$}
\begin{figure}[H]
\includegraphics[width=8cm,height=5cm]{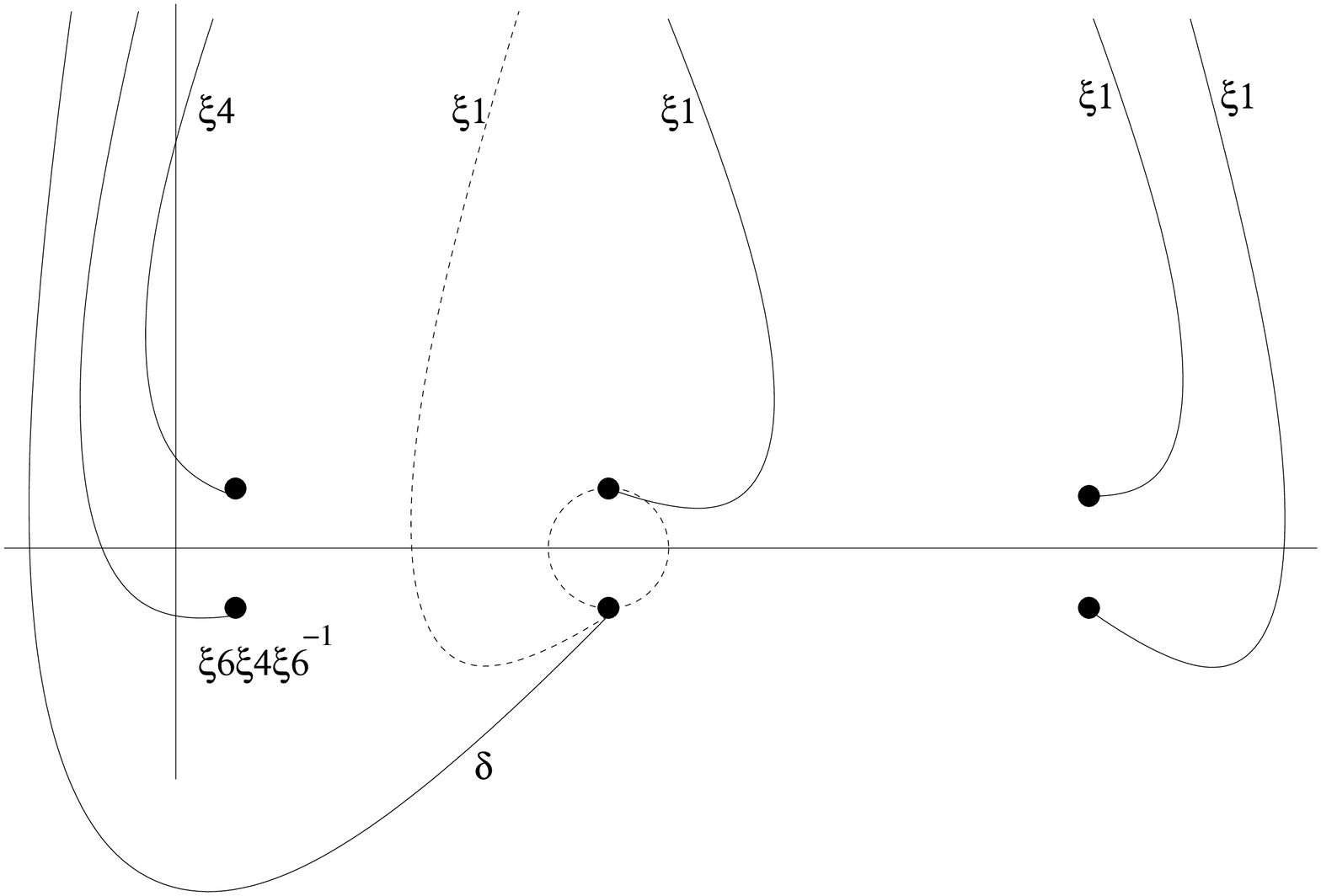}
\caption{\label{figa92a4eta6bis}Generators at $x=\eta_6+\varepsilon$}
\end{figure}

We recall that, at $x=\eta_6-\varepsilon$, the generators are as in Figure~\ref{figa92a4eta6}.
When $x$ makes half--turn counter--clockwise on the circle $\mathbb{S}_\varepsilon(\eta_6)$, from $\eta_6-\varepsilon$ to $\eta_6+\varepsilon$, the $\xi_k$'s are deformed as shown in Figure~\ref{figa92a4eta6bis}, where 
\begin{equation*}
\delta:=(\xi_6\xi_4\xi_6^{-1}\xi_4)\cdot\xi_1\cdot (\xi_6\xi_4\xi_6^{-1}\xi_4)^{-1}.
\end{equation*}

\begin{figure}[H]
\includegraphics[width=8cm,height=5cm]{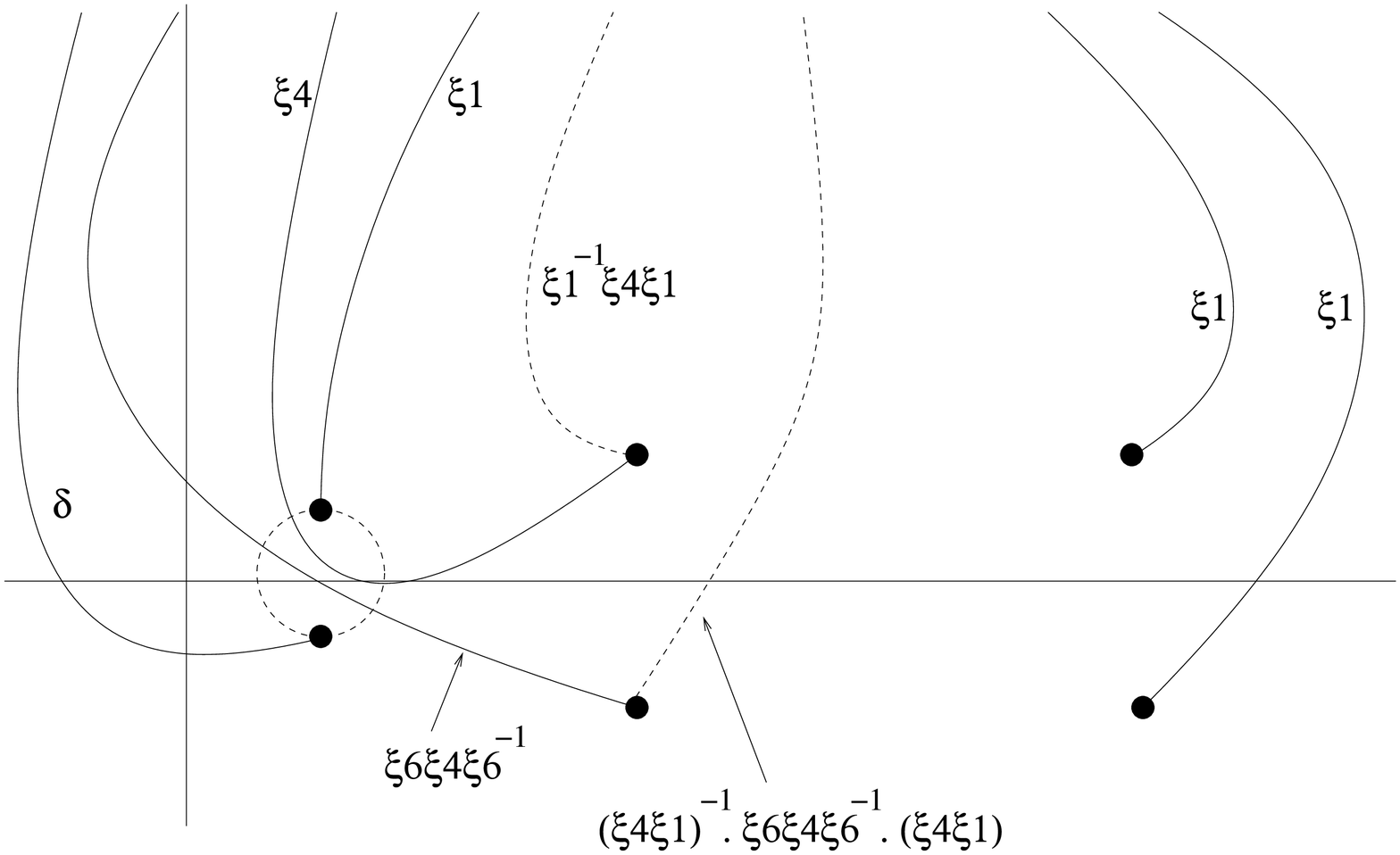}
\caption{\label{figa92a4eta7}Generators at $x=\eta_7-\varepsilon$}
\end{figure}

\begin{lemma}\label{lemmaa92a4}
When $x$ moves on the real axis from $\eta_6+\varepsilon$ to $\eta_7-\varepsilon$, the $\xi_k$'s are deformed as Figure~\ref{figa92a4eta7}.
\end{lemma}

Lemma \ref{lemmaa92a4} is not obvious and will be proved at the end of this section. Before, let us complete the calculation of $\pi_1(\mathbb{CP}^2\setminus C)$.

The monodromy relation at $x=\eta_7$ is a multiplicity $2$ tangent relation given~by
\begin{equation}
\delta=\xi_1.
\end{equation}
Notice that, by (\ref{rela92a4s1}), this relation can also be written as 
\begin{equation}\label{rela92a4eta7s}
\xi_4\xi_1\xi_4\xi_1 = \xi_1\xi_4\xi_1\xi_4.
\end{equation}
Moreover, still using (\ref{rela92a4s1}), we can write the relation (\ref{rela92a4eta400}) as
\begin{displaymath}
\xi_4=(\xi_4\xi_1)^5\cdot \xi_4\cdot (\xi_4\xi_1)^{-5}.
\end{displaymath}
The latter, combined with (\ref{rela92a4eta7s}), implies
\begin{equation*}
\xi_4\xi_1=\xi_1\xi_4.
\end{equation*}
This already shows that the fundamental group $\pi_1(\mathbb{CP}^2\setminus C)$ is abelian. (We do not need to consider the monodromy relations around $L_{\eta_8}$, $L_{\eta_9}$ and $L_{\eta_{10}}$.)

To complete the calculation it remains to prove Lemma \ref{lemmaa92a4}.

\begin{proof}[Proof of Lemma \ref{lemmaa92a4}]
We consider the polynomial
\begin{displaymath}
h(x,u,v):=f(x,u+i\, v)
\end{displaymath}
for $x$, $u$ and $v$ real. We denote by $f_e(x,u,v)$ and $f_o(x,u,v)$ the real and the imaginary part of $h(x,u,v)$ respectively. They have degree $6$ and $5$ respectively in $v$. Suppose there exists a number $x_0\in [\eta_6+\varepsilon,\eta_7-\varepsilon]$ such that $4$ complex solutions of the equation $f(x_0,y)=0$ are on the same vertical line $u=u_0$ in the complex plane $(\mathbb{C},y=u+i\, v)$. This implies that the equations
\begin{displaymath}
f_e(x_0,u_0,v)=f_o(x_0,u_0,v)=0
\end{displaymath}
have $4$ common real solutions $v_1$, $v_2$, $v_3$ and $v_4$. These solutions are non--zero since the equation $\Delta_y(f)(x)=0$ does not have any solution in $[\eta_6+\varepsilon,\eta_7-\varepsilon]$. Therefore, the equations
\begin{displaymath}
f_e(x_0,u_0,v)=f_{oo}(x_0,u_0,v)=0,
\end{displaymath}
where $f_{oo}(x,u,v)=f_o(x,u,v)/v$ (notice that $v$ divides $f_o(x,u,v)$, and $f_{oo}(x,u,v)$ is then a polynomial), also have $v_1$, $v_2$, $v_3$ and $v_4$ as common solutions. As $f_{oo}$ has degree $4$ in $v$, it follows that $f_{oo}(x_0,u_0,v)$ divides $f_e(x_0,u_0,v)$. Therefore the remainder $R(x,u,v)$ of $f_e$ by $f_{oo}$, as a polynomial in $v$, must be identically zero for $u=u_0$ and $x=x_0$ (of course, $R$ is written as $R=R'/R''$ where $R'$ is a polynomial in $x$, $u$ and $v$ while $R''$ is a polynomial depending just on $x$ and $u$). By an easy computation, we see that $R=(R_2'/R_2'')\, v^2 + (R_0'/R_0'')$, where $R_2'$, $R_2''$, $R_0'$ and $R_0''$ are polynomials in $x$ and $u$. Thus, $(x_0,u_0)$ is a common real solution of the equations 
\begin{equation}\label{R}
R_2'(x,u)=R_0'(x,u)=0.
\end{equation}
This implies that $x_0$ is a root of the resultant $\mbox{Res}_u(R_2',R_0')$ of $R_2'$ and $R_0'$ as polynomials in $u$. Note that the condition $\mbox{Res}_u(R_2',R_0')(x_0)=0$ is necessary to have a real partner $u_0$ such that $R_2'(x_0,u_0)=R_0'(x_0,u_0)=0$ but it is not sufficient since the possible partner $u_0$ might be non--real. There are $5$ real solutions $x_{01}, \ldots, x_{05}$ of the equation $\mbox{Res}_u(R_2',R_0')(x)=0$ in the interval $[\eta_6+\varepsilon,\eta_7-\varepsilon]$. Each of them gives a real number, say $u_{0j}$ ($1\leq j\leq 5$), such that $(x_{0j},u_{0j})$ is a solution of (\ref{R}). Now, we have to check if these $5$ solutions give $4$ real roots of the polynomial $v\mapsto f_{oo}(x_0,u_0,v)$. Only the solution $(x_0,u_0):\approx(0.1205,0.0075)$ satisfies this requirement. Therefore we can have one (and only one) overcrossing of $4$ complex roots. To check if it is the case, we look at the solutions $y$ of the equation, in $y$, $f(x,y)=0$ for values of $x$ near $x_0$. \texttt{Maple} actually gives an overcrossing.
\end{proof}

\section{Fundamental group of $\mathbb{CP}^2\setminus C'$}\label{nonabelian}

In this section, we prove that $\pi_1(\mathbb{CP}^2\setminus C') \simeq \mathbb{D}_{10}\times (\mathbb{Z}/3\mathbb{Z})$. As above, we use Zariski--van Kampen's theorem with the pencil given by the vertical lines $L_{\eta}\colon x=\eta$, $\eta\in\mathbb{C}$. 

The discriminant $\Delta_y(f)$ of $f$ as a polynomial in $y$ is the following polynomial in $x$:
\begin{equation*}
\Delta_y(f)(x) = a_0\,{x}^{12} \left( {x}^{4}+2\,{x}^{3}-17\,{x}^{2} +18\,x+9 \right) \left( 25\,{x}^{2}-15\,x-9 \right)^{2} \left( x-1 \right)^{10},
\end{equation*}
where $a_0\in\mathbb{Q}\setminus \{0\}$.
This polynomial has exactly 8 distinct complex roots:\smallskip

$\eta_1 \approx -5.5758$,
\quad $\eta_2 \approx -0.3708$,
\quad $\eta_3 \approx - 0.3677$,

$\eta_4 = 0$,
\quad $\eta_5 \approx 0.9708$,
\quad $\eta_6 = 1$,

$\eta_7 \approx 1.9718 - i\, 0.7077$,
\quad $\eta_8 = \bar \eta_7 \approx 1.9718 + i\, 0.7077$.

\smallskip

\noindent
The singular lines of the pencil are the lines
$L_{\eta_j}$ ($1\leq j\leq 8$) corresponding to these 8 roots. The
lines $L_{\eta_4}$ and $L_{\eta_6}$ intersect the curve at its singular points, while all the other singular lines are tangent to the curve. See Figure~\ref{figa92a4narps}.

\begin{figure}[H]
\includegraphics[width=8cm,height=5cm]{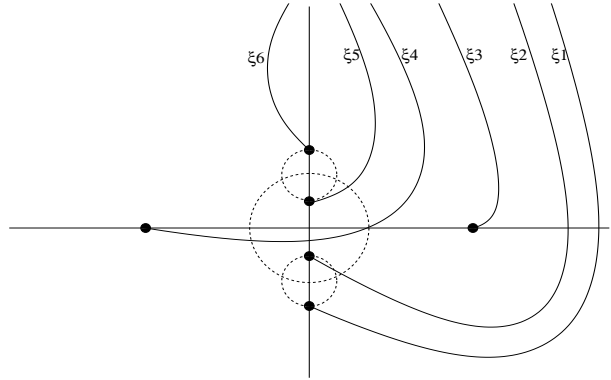}
\caption{\label{figa92a4naeta4}Generators at $x=\eta_4-\varepsilon$}
\end{figure}

Here we start with the generic line $L_{\eta_4-\varepsilon}$ and we
choose generators $\xi_1,\ldots,\xi_6$ of 
$\pi_1(L_{\eta_4-\varepsilon}\setminus C')$ as in Figure~\ref{figa92a4naeta4}, where $\varepsilon>0$ is small enough. (To determine the position of the complex roots of the equation $f(\eta_4-\varepsilon,y)=0$ one may use (\ref{b0}) and (\ref{pp0}) below.)

As above, to find the monodromy relations around the singular lines $L_{\eta_j}$ ($1\leq j\leq 8$) of the pencil, we fix a `standard' system of generators $\sigma_1, \ldots,\sigma_8$ of the fundamental group $\pi_1(\mathbb{C}\setminus \{\eta_1,\ldots, \eta_8\})$, where each $\sigma_j$ is a lasso (oriented counter--clockwise) around $\eta_j$ and with base point $\eta_4-\varepsilon$. For $j\not= 7,8$, the tail of $\sigma_j$ is a union of real segments and half--circles around the exceptional parameters $\eta_l$ ($l\not=j$) located in the real axis between the base point $\eta_4-\varepsilon$ and $\eta_j$, while its head is the circle $\mathbb{S}_\varepsilon(\eta_j)$. The lasso $\sigma_7$ corresponding to the non--real root $\eta_7$ has the form $\zeta\theta\zeta^{-1}$, where $\theta$ is the loop obtained by moving $x$ once on the circle $\mathbb{S}_\varepsilon(\eta_7)$, starting at $\Re(\eta_7)+i\, (\Im(\eta_7)+\varepsilon)$, and $\zeta$ the path obtained when $x$ makes half--turn on the circle $\mathbb{S}_\varepsilon(\eta_4)$, from $\eta_4-\varepsilon$ to $\eta_4+\varepsilon$, moves on the real axis from $\eta_4+\varepsilon$ to $\eta_5-\varepsilon$, makes half--turn on the circle $\mathbb{S}_\varepsilon(\eta_5)$, from $\eta_5-\varepsilon$ to $\eta_5+\varepsilon$, moves on the real axis from $\eta_5+\varepsilon$ to $\eta_6-\varepsilon$, makes half--turn on the circle $\mathbb{S}_\varepsilon(\eta_6)$, from $\eta_6-\varepsilon$ to $\eta_6+\varepsilon$, moves on the real axis from $\eta_6+\varepsilon$ to $\Re(\eta_7)$, and finally moves in a straight line from $\Re(\eta_7)$ to $\Re(\eta_7)+i\, (\Im(\eta_7)+\varepsilon)$. The lasso $\sigma_8$ is defined similarly from a loop $\theta$ and a path $\zeta$ meeting at $\Re(\eta_8)+i\, (\Im(\eta_8)-\varepsilon)$.

The monodromy relations are now given as follows.

\subsection*{Monodromy relations at $x=\eta_3$}
\begin{figure}[H]
\includegraphics[width=6cm,height=6cm]{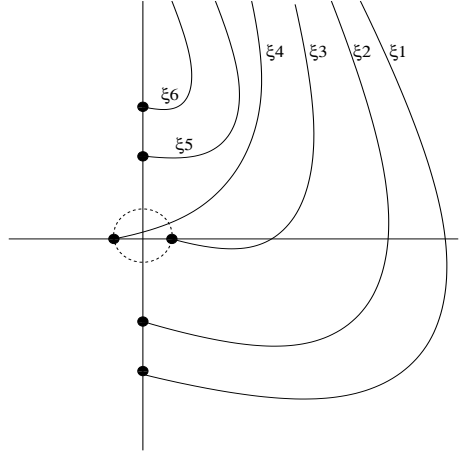}
\caption{\label{figa92a4naeta3}Generators at $x=\eta_3+\varepsilon$}
\end{figure}

When $x$ moves on the real axis from $\eta_4-\varepsilon$ to $\eta_3+\varepsilon$, the $\xi_k$'s are deformed as in Figure~\ref{figa92a4naeta3}. The singular line $L_{\eta_3}$  is tangent to $C'$ at one simple point, with intersection multiplicity 2, so the monodromy relation around this line is simply given by 
\begin{equation}\label{rela92a4naeta3}
\xi_4=\xi_3.
\end{equation}

\subsection*{Monodromy relations at $x=\eta_2$}
\begin{figure}[H]
\includegraphics[width=6cm,height=6cm]{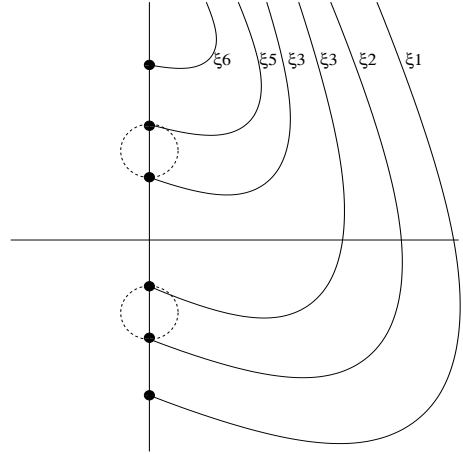}
\caption{\label{figa92a4naeta2}Generators at $x=\eta_2+\varepsilon$}
\end{figure}

Now, when $x$ makes half--turn counter--clockwise on the circle $\mathbb{S}_\varepsilon(\eta_3)$, from $\eta_3+\varepsilon$ to $\eta_3-\varepsilon$, and then moves on the real axis from $\eta_3-\varepsilon$ to $\eta_2+\varepsilon$, the $\xi_k$'s are deformed as in Figure~\ref{figa92a4naeta2}. The singular line $L_{\eta_2}$ is tangent to the curve at two simple points, in both case with intersection multiplicity 2. The monodromy relations around this line are given by 
\begin{eqnarray}
\label{rela92a4naeta21} && \xi_5=\xi_3,\\
\label{rela92a4naeta22} && \xi_3=\xi_2.
\end{eqnarray}

\begin{figure}[H]
\includegraphics[width=8cm,height=5cm]{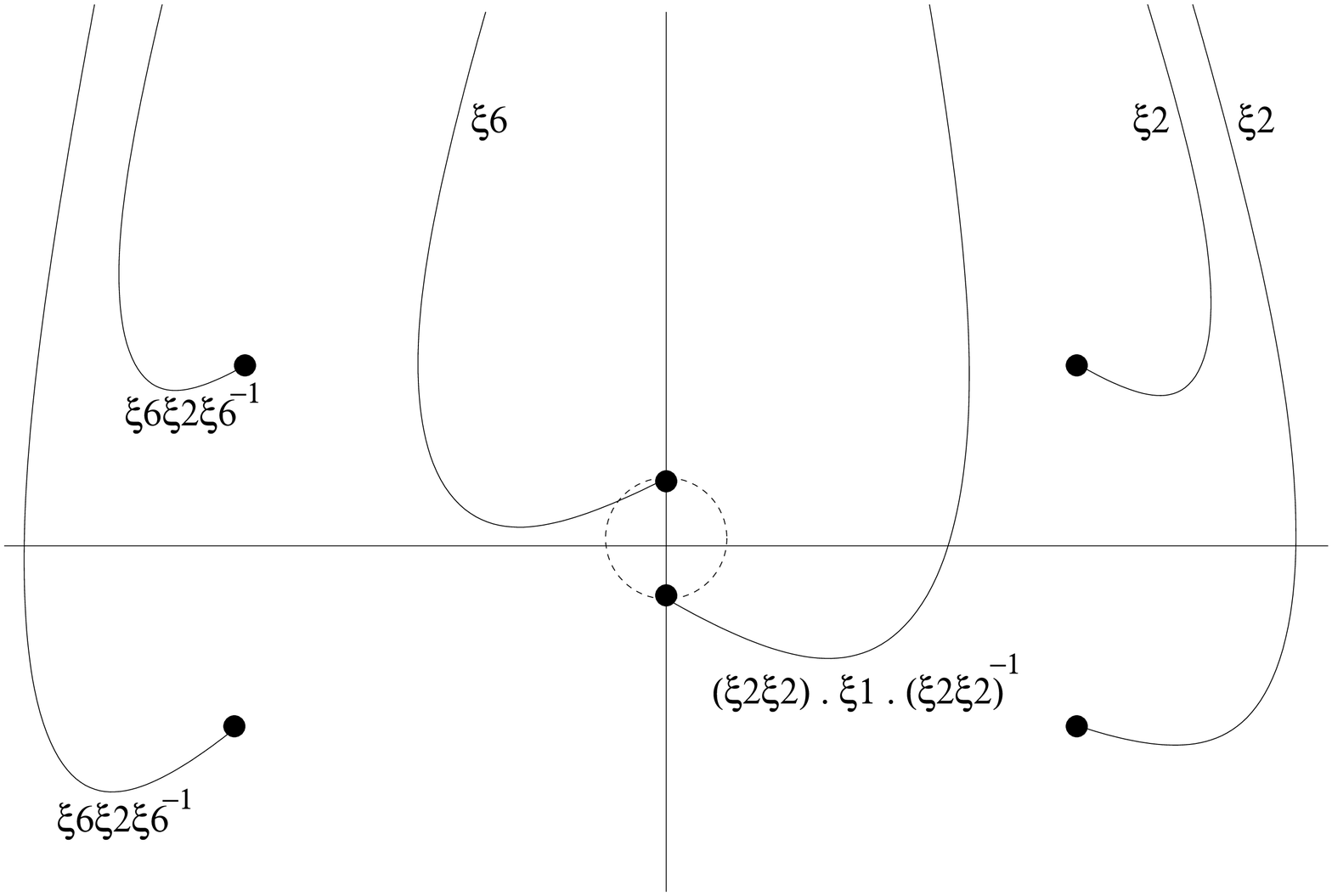}
\caption{\label{figa92a4naeta1}Generators at $x=\eta_1+\varepsilon$}
\end{figure}

\subsection*{Monodromy relations at $x=\eta_1$}
In Figure~\ref{figa92a4naeta1}, we show how the $\xi_k$'s are deformed when $x$ makes half--turn counter--clockwise on the circle $\mathbb{S}_\varepsilon(\eta_2)$, from $\eta_2+\varepsilon$ to $\eta_2-\varepsilon$, and then moves on the real axis from $\eta_2-\varepsilon$ to $\eta_1+\varepsilon$. The monodromy relation around $L_{\eta_1}$ is a multiplicity 2 tangent relation: 
\begin{equation}\label{rela92a4naeta1}
\xi_6=\xi_2^2\xi_1\xi_2^{-2}.
\end{equation}

Notice that the relations (\ref{rela92a4naeta3})--(\ref{rela92a4naeta1}) show that the \emph{vanishing relation at infinity} can be written as
\begin{equation}\label{rela92a4nainfini}
(\xi_2\xi_2\xi_1)^2=e,
\end{equation}
where $e$ is the unit element.

\subsection*{Monodromy relations at $x=\eta_4$}
\begin{figure}[H]
\includegraphics[width=8cm,height=5cm]{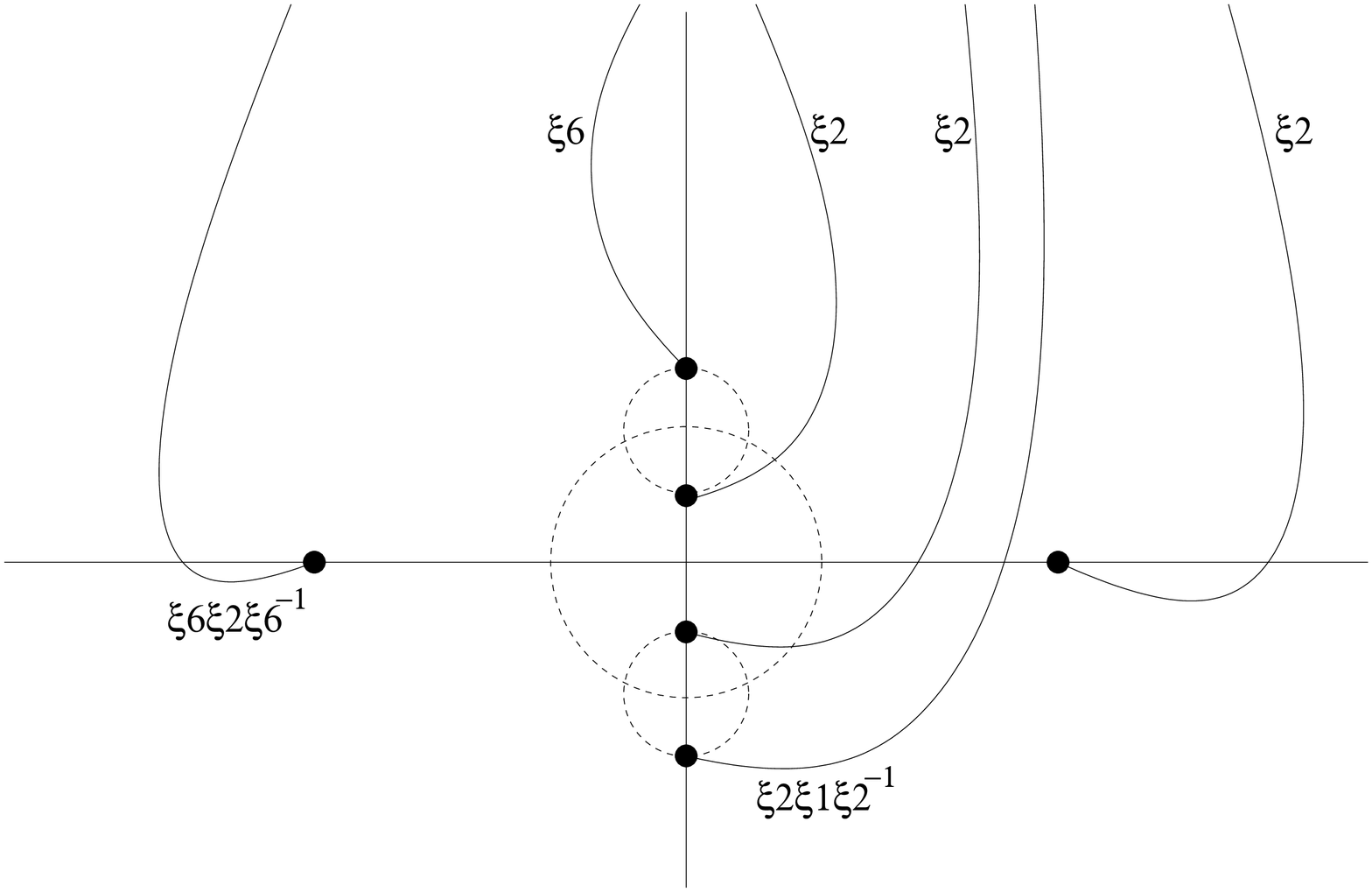}
\caption{\label{figa92a4naeta4bis}Generators at $x=\eta_4-\varepsilon$}
\end{figure}

At $x=\eta_4-\varepsilon$, the generators are shown in Figure~\ref{figa92a4naeta4}.
By (\ref{rela92a4naeta3}), (\ref{rela92a4naeta21}) and (\ref{rela92a4naeta22}), Figure~\ref{figa92a4naeta4} is the same as Figure~\ref{figa92a4naeta4bis}. The singular line $L_{\eta_4}$ passes through the origin which a type $\textbf{A}_9$ singular point of the curve. At this point the curve has two branches $K_+$ and $K_-$ given by
\begin{equation}\label{b0}
K_{\pm}:\quad x=y^2-{\frac{1}{3}}\, y^4\pm{\frac{2}{9}}\, i\, \sqrt{3}\, y^5+\, \mbox{higher terms.}
\end{equation}
An easy computation shows that Puiseux parametrizations of these branches are:
\begin{equation}\label{pp0}
K_{\pm}:\quad x=t^2,\quad y=t+{\frac{1}{6}}\, t^3\mp{\frac{\sqrt{3}}{9}}\, i\, t^4+\, \mbox{higher terms.}
\end{equation}
The monodromy relations at $x=\eta_4$ are then given by
\begin{eqnarray}
\label{rela92a4naeta41} & & \xi_2\xi_1\xi_2^{-1} = (\xi_6\xi_2)\cdot \xi_6\xi_2\xi_6^{-1}\cdot (\xi_6\xi_2)^{-1},\\ 
&\label{rela92a4naeta42}  & \xi_2 = (\xi_6\xi_2)^2\cdot\xi_6\cdot (\xi_6\xi_2)^{-2},\\
&\label{rela92a4naeta43}  & \xi_2=(\xi_6\xi_2^3\xi_1)\cdot \xi_2\xi_1\xi_2^{-1}\cdot (\xi_6\xi_2^3\xi_1)^{-1},\\
&\label{rela92a4naeta44}  & \xi_6=(\xi_6\xi_2^3\xi_1)\cdot (\xi_2\xi_1)\cdot\xi_2\cdot(\xi_2\xi_1)^{-1}\cdot (\xi_6\xi_2^3\xi_1)^{-1}.
\end{eqnarray}
Note that, by (\ref{rela92a4naeta1}), all of them are equivalent to
\begin{eqnarray}\label{rela92a4nas1}
\xi_1\xi_2\xi_1\xi_2\xi_1 = \xi_2\xi_1\xi_2\xi_1\xi_2.
\end{eqnarray}

\subsection*{Monodromy relations at $x=\eta_5$}
\begin{figure}[H]
\includegraphics[width=8cm,height=5cm]{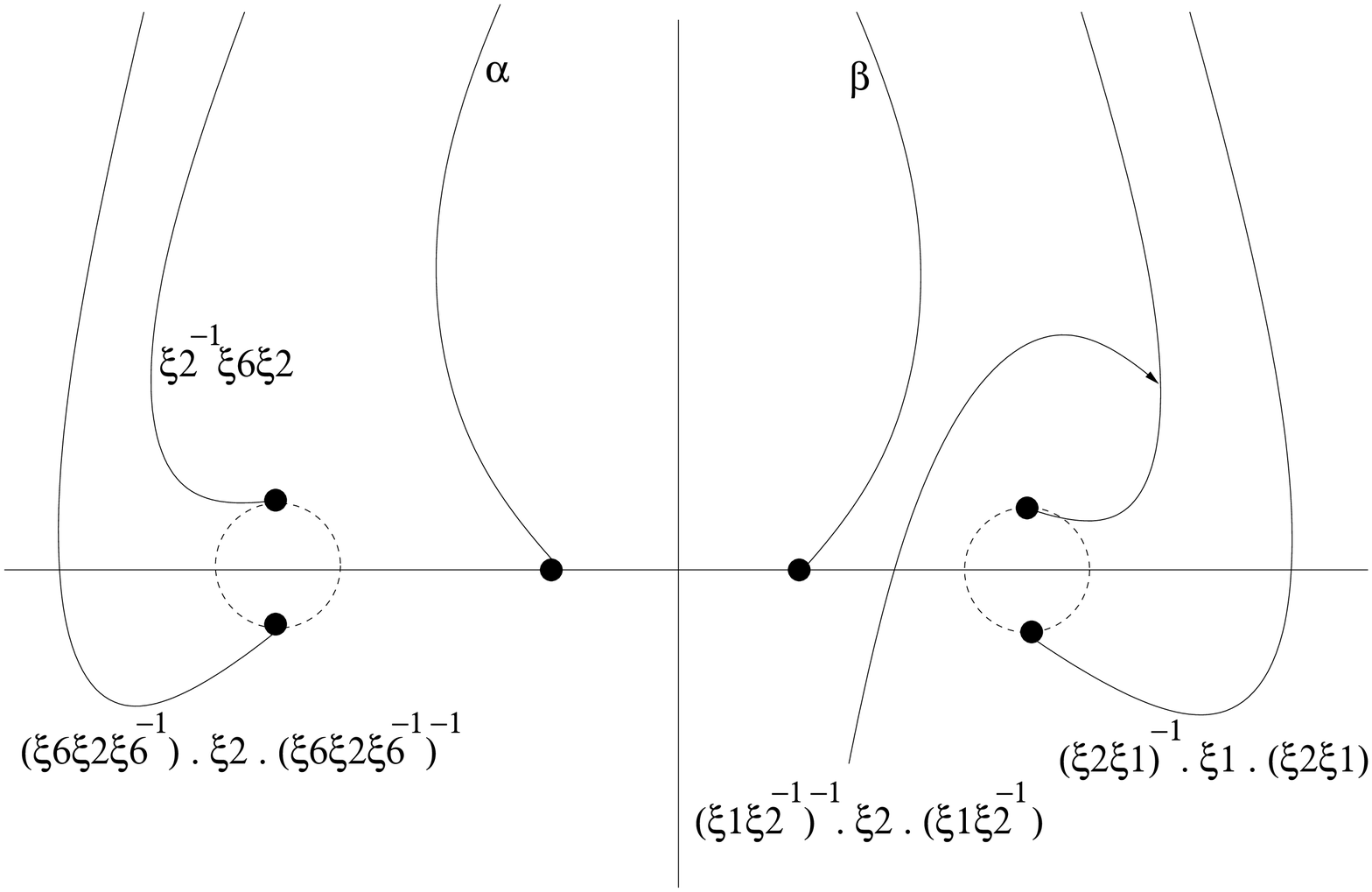}
\caption{\label{figa92a4naeta5}Generators at $x=\eta_5-\varepsilon$}
\end{figure}

At $x=\eta_4-\varepsilon$, the generators are as in Figure~\ref{figa92a4naeta4bis}.
Now, when $x$ makes half--turn counter--clockwise on the circle $\mathbb{S}_\varepsilon(\eta_4)$, from $\eta_4-\varepsilon$ to $\eta_4+\varepsilon$, then moves on the real axis from $\eta_4+\varepsilon$ to $\eta_5-\varepsilon$, the $\xi_k$'s are deformed as in Figure~\ref{figa92a4naeta5}, where
\begin{eqnarray*}
&& \alpha:=(\xi_2^{-1}\xi_6\xi_2)^{-1}\cdot \xi_6\xi_2\xi_6^{-1}\cdot (\xi_2^{-1}\xi_6\xi_2), \\
&& \beta:=(\xi_1\xi_2^{-1})^{-1}\cdot (\xi_2\xi_1)\cdot \xi_2\cdot (\xi_2\xi_1)^{-1}
\cdot (\xi_1\xi_2^{-1}).
\end{eqnarray*}
The singular line $L_{\eta_1}$ is tangent to the curve at two simple points, in both cases with intersection multiplicity 2. The monodromy relations at $x=\eta_5$ are then given as follows:
\begin{eqnarray}
\label{rela92a4naeta51} & & (\xi_1\xi_2^{-1})^{-1}\cdot \xi_2\cdot (\xi_1\xi_2^{-1}) = (\xi_2\xi_1)^{-1}\cdot \xi_1\cdot (\xi_2\xi_1),\\ 
&\label{rela92a4naeta52}  & (\xi_6\xi_2\xi_6^{-1})\cdot \xi_2\cdot (\xi_6\xi_2\xi_6^{-1})^{-1} = \xi_2^{-1}\xi_6\xi_2.
\end{eqnarray}
By (\ref{rela92a4naeta1}) and (\ref{rela92a4nas1}), both of them are equivalent to 
\begin{eqnarray}\label{rela92a4nas2}
\xi_2\xi_1\xi_1 = \xi_1\xi_1\xi_2.
\end{eqnarray}

\subsection*{Monodromy relations at $x=\eta_6$}
\begin{figure}[H]
\includegraphics[width=8cm,height=5cm]{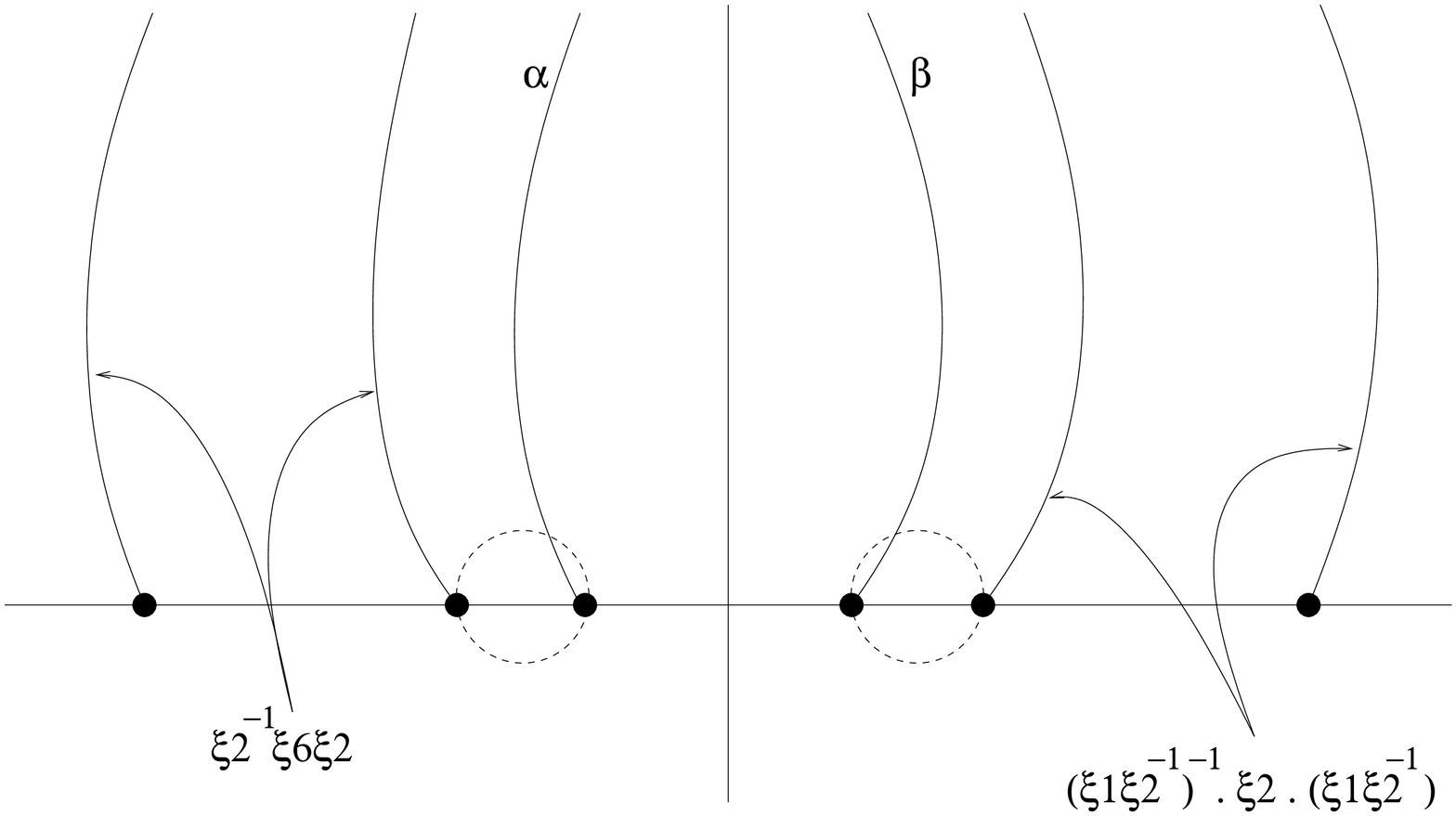}
\caption{\label{figa92a4naeta6}Generators at $x=\eta_6-\varepsilon$}
\end{figure}

In Figure~\ref{figa92a4naeta6}, we show how the generators are deformed when $x$ makes half--turn counter--clockwise on the circle $\mathbb{S}_\varepsilon(\eta_5)$, from $\eta_5-\varepsilon$ to $\eta_5+\varepsilon$, and then moves on the real axis from $\eta_5+\varepsilon$ to $\eta_6-\varepsilon$. The singular line~$L_{\eta_6}$ passes through the points $(1,1)$ and $(1,-1)$ which are both singularities of type $\textbf{A}_4$. Puiseux parametrizations of the curve at these points are given by:
\begin{eqnarray*}
(1,1):\quad & x=1+t^2,\ & y=1+t^4+3\, i\, \sqrt{3}\, t^5 + \mbox{higher terms},\\
(1,-1):\quad & x=1+t^2,\ & y=-1-t^4-3\, i\, \sqrt{3}\, t^5 + \mbox{higher terms}.
\end{eqnarray*}
The monodromy relations at $x=\eta_6$ are then given by
\begin{eqnarray}
\label{rela92a4naeta61} & & \xi_2\cdot \xi_1\xi_2\xi_1^{-1}\cdot \xi_2\cdot \xi_1\xi_2\xi_1^{-1}\cdot \xi_2 = \xi_1\xi_2\xi_1^{-1}\cdot \xi_2\cdot \xi_1\xi_2\xi_1^{-1}\cdot \xi_2\cdot\xi_1\xi_2\xi_1^{-1},\\
&\label{rela92a4naeta62}  & \alpha = (\xi_2^{-1}\xi_6\xi_2\alpha)^2 \cdot \xi_2^{-1}\xi_6\xi_2\cdot (\xi_2^{-1}\xi_6\xi_2\alpha)^{-2}.
\end{eqnarray}
Notice that (\ref{rela92a4naeta61}) is automatically satisfied, while  (\ref{rela92a4naeta62}) can be written as 
\begin{eqnarray}\label{rela92a4nas3}
(\xi_2\xi_1)^4 = \xi_1\xi_2
\end{eqnarray}
or, equivalently, as 
\begin{eqnarray}\label{rela92a4nas33}
(\xi_2\xi_1\xi_1)^2 = e.
\end{eqnarray}
Indeed, using (\ref{rela92a4nas2}) under the form $\xi_1^{-1}\xi_2\xi_1 = \xi_1\xi_2\xi_1^{-1}$, the relation (\ref{rela92a4naeta61}) can be written as 
\begin{eqnarray*}
\xi_2 \cdot (\xi_1^{-1}\xi_2\xi_1 \cdot \xi_2 \cdot \xi_1)\xi_2\xi_1^{-1}  \cdot \xi_2 = \xi_1^{-1}\xi_2\xi_1 \cdot \xi_2 \cdot (\xi_1^{-1}\xi_2\xi_1 \cdot \xi_2 \cdot \xi_1)\xi_2\xi_1^{-1}, 
\end{eqnarray*}
which is nothing but 
$\xi_1\xi_2\cdot \xi_2\xi_1\xi_2 = \xi_2\xi_1\cdot \xi_2\xi_2\xi_1$ by (\ref{rela92a4nas1}).
By (\ref{rela92a4nainfini}) this equality is always satisfied.
The relation (\ref{rela92a4naeta62}) is written as 
\begin{eqnarray*}
\xi_6\xi_2\xi_6^{-1}\cdot (\xi_2^{-1}\xi_6\xi_2\cdot \xi_6\xi_2\xi_6^{-1})^2 =
(\xi_2^{-1}\xi_6\xi_2\cdot \xi_6\xi_2\xi_6^{-1})^2 \cdot \xi_2^{-1}\xi_6\xi_2,
\end{eqnarray*}
which is the same as
\begin{eqnarray*}
\xi_6\xi_2\xi_6^{-1}\cdot (\xi_6\xi_2\xi_6^{-1}\cdot \xi_2)^2 =
(\xi_6\xi_2\xi_6^{-1}\cdot \xi_2)^2 \cdot \xi_2^{-1}\xi_6\xi_2
\end{eqnarray*}
by (\ref{rela92a4naeta52}). Equivalently
$\xi_6\xi_2\xi_6^{-1}\cdot \xi_2 \cdot \xi_6\xi_2\xi_6^{-1} = \xi_2\xi_6\xi_2$.
By (\ref{rela92a4naeta1}), $\xi_6$ can be eliminated so 
\begin{eqnarray}\label{rela92a4nas5}
\xi_1\xi_2\cdot (\xi_2\xi_1)^{-1}\cdot \xi_2\xi_2\xi_1 = \xi_2\xi_1\xi_2\xi_1\xi_2^{-1}.
\end{eqnarray}
By (\ref{rela92a4nainfini}) this is the same as $\xi_1\xi_2 = (\xi_2\xi_1)^4$, while (\ref{rela92a4nainfini}) and (\ref{rela92a4nas1}) show that (\ref{rela92a4nas5}) is also the same as $(\xi_2\xi_1\xi_1)^2=e$.

It remains to find the monodromy relations around the singular lines $L_{\eta_7}$ and $L_{\eta_8}$ corresponding to the non--real roots $\eta_7$ and $\eta_8$ of the discriminant $\Delta_y(f)(x)$.

\subsection*{Monodromy relations at $x=\eta_7$ and $x=\eta_8$}
\begin{figure}[H]
\includegraphics[width=8cm,height=5cm]{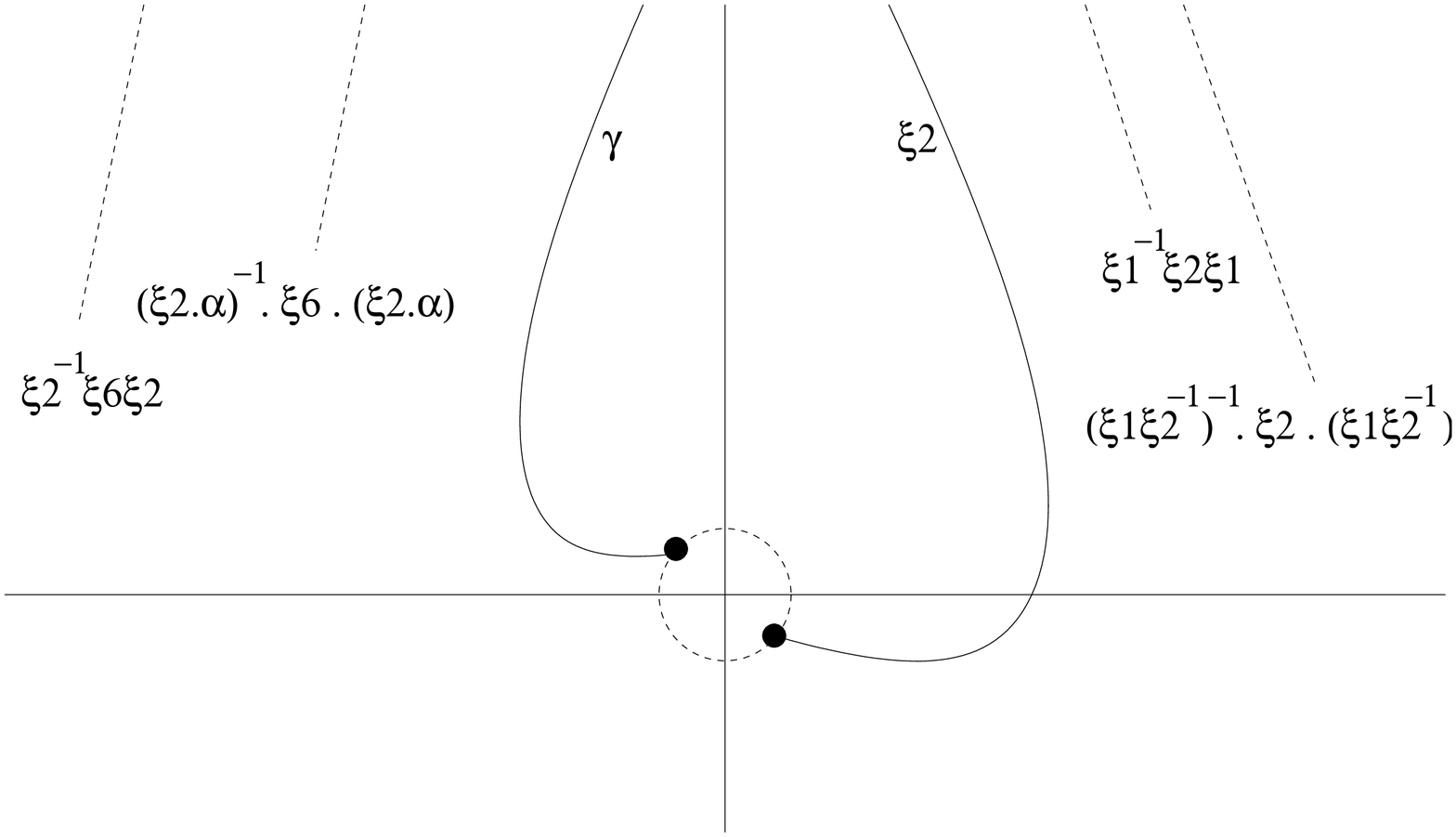}
\caption{\label{figa92a4naeta7}Generators at $x=\Re(\eta_7)+i\, (\Im(\eta_7)+\varepsilon)$}
\end{figure}

\begin{figure}[H]
\includegraphics[width=8cm,height=5cm]{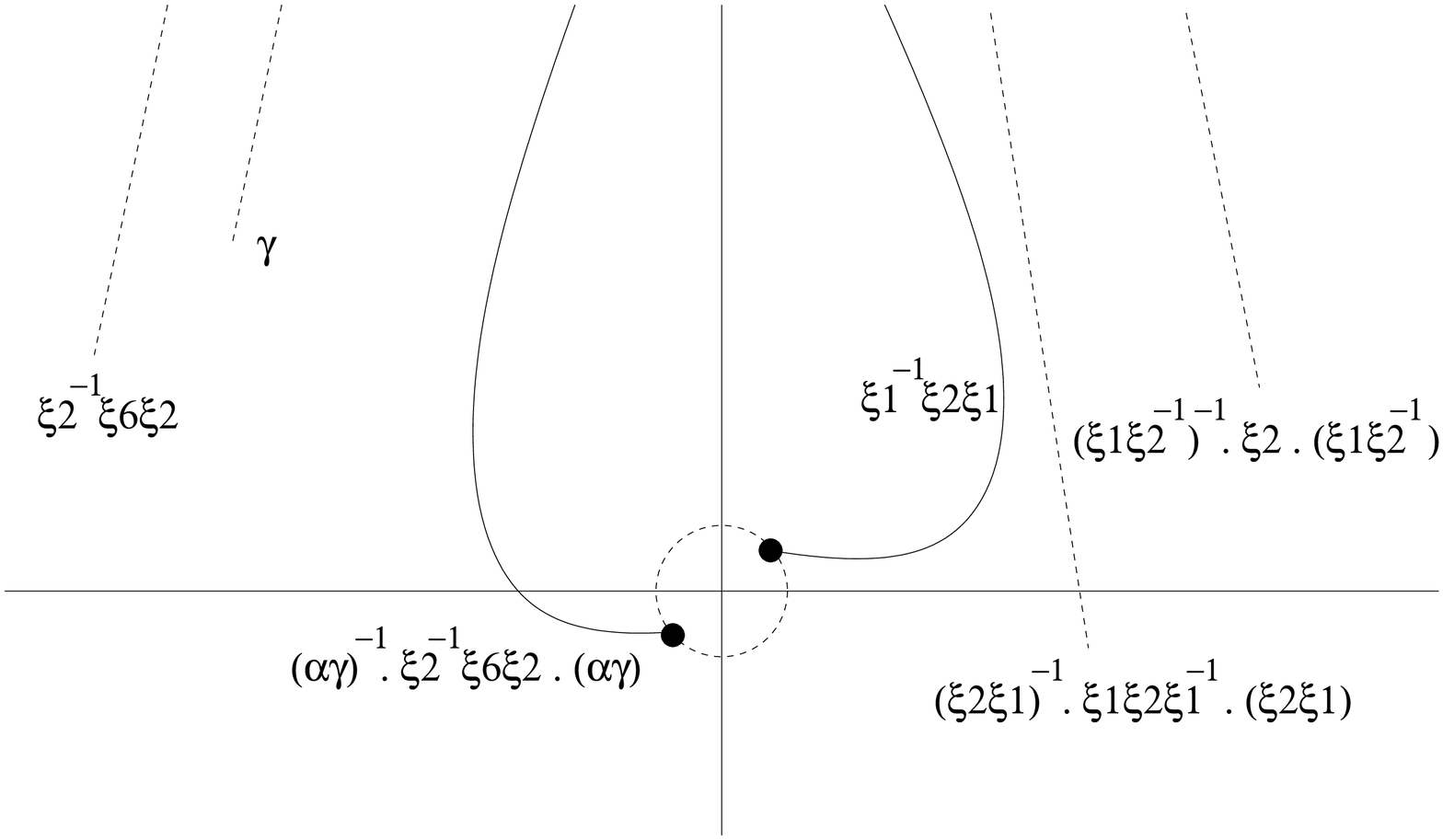}
\caption{\label{figa92a4naeta8}Generators at $x=\Re(\eta_8)+i\, (\Im(\eta_8)-\varepsilon)$}
\end{figure}

Figure~\ref{figa92a4naeta6} shows the generators at $x=\eta_6-\varepsilon$.
In Figure~\ref{figa92a4naeta7} (respectively Figure~\ref{figa92a4naeta8}), we show how the $\xi_k$'s are deformed when $x$ makes half--turn counter--clockwise on the circle $\mathbb{S}_\varepsilon(\eta_6)$, from $\eta_6-\varepsilon$ to $\eta_6+\varepsilon$, then moves on the real axis from $\eta_6+\varepsilon$ to $\Re(\eta_7)$, and finally moves straight along the line $(\Re(\eta_7),\eta_7)$ from $\Re(\eta_7)$ to $\Re(\eta_7)+i\, (\Im(\eta_7)+\varepsilon)$ (respectively along the line $(\Re(\eta_8),\eta_8)$ from $\Re(\eta_8)$ to $\Re(\eta_8)+i\, (\Im(\eta_8)-\varepsilon)$), where 
\begin{equation*}
\gamma:=(\xi_2^{-1}\xi_6\xi_2\alpha)^{-1} \cdot \alpha \cdot (\xi_2^{-1}\xi_6\xi_2\alpha).
\end{equation*} 
(In these figures we concentrate only on the generators which may give \emph{a priori} some relations.) The monodromy relations at $x=\eta_7$ and $x=\eta_8$ are multiplicity 2 tangent relations given by
\begin{equation}
\label{rela92a4naeta7} \gamma=\xi_2,
\end{equation}
and
\begin{equation}
\label{rela92a4naeta8} (\alpha \gamma)^{-1} \cdot  \xi_2^{-1}\xi_6\xi_2 \cdot (\alpha \gamma)=\xi_1^{-1}\xi_2\xi_1,
\end{equation}
respectively.
In fact, these relations are automatically satisfied. Indeed, the relation (\ref{rela92a4naeta7}) is written as
\begin{equation*}
\xi_6\xi_2\xi_6^{-1} \cdot (\xi_2^{-1}\xi_6\xi_2\cdot \xi_6\xi_2\xi_6^{-1}) \cdot \xi_2^{-1}\xi_6 = (\xi_2^{-1}\xi_6\xi_2\cdot \xi_6\xi_2\xi_6^{-1}) \cdot \xi_2^{-1}\xi_6\xi_2.
\end{equation*}
But, by (\ref{rela92a4naeta52}), we know that $\xi_2^{-1}\xi_6\xi_2\cdot \xi_6\xi_2\xi_6^{-1} = \xi_6\xi_2\xi_6^{-1} \cdot \xi_2$, so the relation above is always satisfied. Now, using (\ref{rela92a4naeta7}), the relation (\ref{rela92a4naeta8}) is written as
\begin{equation*}
(\xi_2^{-1}\xi_6\xi_2\cdot \xi_6\xi_2\xi_6^{-1}) \cdot \xi_2^{-1}\xi_6\xi_2 \cdot \xi_2 = \xi_6\xi_2\xi_6^{-1} \cdot \xi_2^{-1}\xi_6\xi_2 \cdot  \xi_2 \cdot \xi_1^{-1}\xi_2\xi_1,
\end{equation*}
which is equivalent to $\xi_2\xi_6\xi_2\xi_2 = \xi_6\xi_2\xi_2\xi_1^{-1}\xi_2\xi_1$ by (\ref{rela92a4naeta52}). The latter is always satisfied, by (\ref{rela92a4naeta1}).

Altogether, we have proved that the fundamental group $\pi_1(\mathbb{CP}^2\setminus C')$ is presented by the generators $\xi_1$ and $\xi_2$ and the relations (\ref{rela92a4nainfini}), (\ref{rela92a4nas1}), (\ref{rela92a4nas2}) and (\ref{rela92a4nas3}). The relation (\ref{rela92a4nainfini}) can be written as
\begin{equation}\label{rela92a4nas4}
(\xi_2\xi_1\xi_2)^2=e.
\end{equation}
This shows that (\ref{rela92a4nas1}) is equivalent to (\ref{rela92a4nas3}). The relation (\ref{rela92a4nas2}) is automatically satisfied. Indeed, by (\ref{rela92a4nas3}), it is equivalent to 
\begin{equation*}
(\xi_2\xi_1\xi_1)^2=(\xi_1\xi_2\xi_1)^2.
\end{equation*}
But both sides are equal to $e$, by (\ref{rela92a4nas3}) under the form (\ref{rela92a4nas33}). Hence, $\pi_1(\mathbb{CP}^2\setminus C')$ is simply presented by the generators $\xi_1$ and $\xi_2$ and the relations (\ref{rela92a4nas3}) and (\ref{rela92a4nas4}). After the change $a:=\xi_2\xi_1\xi_2$ and $b:=\xi_2\xi_1$, the presentation is also given by
\begin{equation*}
\pi_1(\mathbb{CP}^2\setminus C') \simeq \bigl\langle\, a,\, b\, \mid a^2=e,\, aba=b^4\, \bigr\rangle.
\end{equation*}

\begin{lemma}
The generator $b$ satisfies the following two properties:
\begin{enumerate}
\item $b^{15}=e$;
\item $b^5$ is in the centre of
$\pi_1(\mathbb{CP}^2\setminus C')$.
\end{enumerate}
\end{lemma}

\begin{proof}
Since $a^2=e$, the
relation $aba=b^4$ gives $b^{16}=ab^4a=b$,
that is, $b^{15}=e$
as desired. To show that $b^5$ is in the centre of
$\pi_1(\mathbb{CP}^2\setminus C')$ we write:
\begin{eqnarray*}
b^5ab^{-5}a^{-1}  = b\cdot b^4\cdot ab^{-5}a^{-1}
= b\cdot aba\cdot ab^{-5}a^{-1}=\\
 ba\cdot b^{-4}\cdot a^{-1}
= ba\cdot a^{-1}b^{-1}a^{-1}\cdot a^{-1}
= e.
\end{eqnarray*}
\end{proof}

It follows from the lemma that $\pi_1(\mathbb{CP}^2\setminus C')$ is also presented as:
\begin{eqnarray*}
\pi_1(\mathbb{CP}^2\setminus C') & \simeq
& \bigl\langle 
a,\, b\mid a^2=e,\ aba=b^4,\, b^{15}=e,\, b^5a=ab^5
\bigr\rangle \\
& \simeq & \bigl\langle
a,\, b,\, c,\, d\mid a^2=b^{15}=e,\, aba=b^4,\, b^5a=ab^5,\, c=b^6, \\
&& \qquad d=b^5,\,  da=ad,\, db=bd,\, dc=cd
\bigr\rangle \\
& \simeq & \bigl\langle
a,\, b,\, c,\, d\mid a^2=b^{15}=e,\, aba=b^4,\, c=b^6,\, d=b^5,\\
&& \qquad b=cd^{-1},\,  da=ad,\, db=bd,\, dc=cd
\bigr\rangle \\
& \simeq & \bigl\langle
a,\, c,\, d\mid a^2=c^5=d^3=e,\, acd^{-1}a=c^4d^{-1},\,  da=ad,\, dc=cd
\bigr\rangle \\
& \simeq & \bigl\langle
a,\, c,\, d\mid a^2=c^5=d^3=e,\, aca=c^4,\,  da=ad,\, dc=cd
\bigr\rangle \\
& \simeq & \mathbb{D}_{10}\times (\mathbb{Z}/3\mathbb{Z}).
\end{eqnarray*}

\bibliographystyle{amsplain}

\end{document}